\documentclass[12pt, twoside]{article}
\usepackage{amsmath,amsthm}
\usepackage{amssymb,latexsym}
\usepackage{enumerate}
\usepackage[mathscr]{eucal}
\usepackage{bbm}

\newtheorem{thm}{Theorem}

\newtheorem{thms}{Theorem}[section]
\newtheorem{props}[thms]{Proposition}

\newtheorem*{xlem}{Lemma}

\theoremstyle{definition}

\newtheorem{exa}{Example}

\newtheorem*{xrem}{Remark}

\numberwithin{equation}{section}

\multlinetaggap\parindent

\def\ods{\vskip-\lastskip\vskip4pt plus2pt}

\newcommand{\vg}{\vadjust{\goodbreak}}

\newcommand{\textbot}[1]{\begin{minipage}[b]{0.8\linewidth}#1
\end{minipage}}

\newcommand{\calN}{\mathcal N}
\newcommand{\calL}{\mathcal L}
\newcommand{\calB}{\mathcal B}
\newcommand{\calD}{\mathcal D}
\newcommand{\calE}{\mathcal E}
\newcommand{\calG}{\mathcal G}
\newcommand{\calO}{\mathcal O}
\newcommand{\calS}{\mathcal S}

\newcommand{\symC}{\mathbb C}
\newcommand{\Cm}{{\mathbb C}^m}
\newcommand{\Cn}{{\mathbb C}^n}
\newcommand{\symR}{\mathbb R}
\newcommand{\Rn}{{\mathbb R}^n}
\newcommand{\symN}{\mathbb N}

\newcommand{\B}{{\mathbf B}}

\newcommand{\Mmm}{M_{m\times m}}
\newcommand{\mm}{m\times m}

\newcommand{\scrF}{{\mathscr F}}
\newcommand{\adj}{{\rm adj}}
\newcommand{\symjed}{\mathbbm{1}}

\def\<{\langle}
\def\>{\rangle}

\DeclareMathOperator{\convsupp}{conv\,supp}
\DeclareMathOperator{\supp}{supp}
\DeclareMathOperator{\hRe}{Re}
\DeclareMathOperator{\hIm}{Im}
\DeclareMathOperator{\sgn}{sgn}

\newenvironment{aitem}{\begin{list}{}%
{\def\makelabel##1{\rlap{\rm##1}\hss}%
\setlength{\labelsep}{5pt}
\settowidth\labelwidth{\rm(b)}
\setlength{\leftmargin}%{\parindent}\addtolength{\leftmargin}%
{\labelwidth}\addtolength{\leftmargin}{\labelsep}
\setlength{\topsep}{4pt}%
\setlength{\parsep}{0pt}%
\setlength{\itemsep}{0pt}%
}}
{\end{list}}

\renewcommand{\[}{\mathclose[}
\renewcommand{\]}{\mathopen]}

\begin{document}

\baselineskip=17pt

\title{One-parameter convolution semigroups
of rapidly decreasing distributions}
\author{Jan Kisy\'nski\\
Institute of Mathematics, Polish Academy of Sciences\\
\'Sniadeckich 8, 00-956 Warszawa, Poland\\
E-mail: jan.kisynski@gmail.com}

\date{}

\maketitle

\renewcommand{\thefootnote}{}



\begin{abstract}
Let $\Mmm$ denote the set of $\mm$ matrices with complex entries, and let $\calG(\partial_1,\ldots,\partial_n)$ be an $\mm$ matrix whose entries are partial differential operators on $\Rn$ with constant complex coefficients. It is proved that $\calG(\partial_1,\ldots,\partial_n)\otimes \delta$ is the generating distribution of a smooth one-parameter convolution semigroup of $\Mmm$-valued rapidly decreasing distributions on $\Rn$ if and only if
$$\sup_{(\xi_1,\ldots,\xi_n)\in\Rn}\hRe\sigma(\calG(i\xi_1,\ldots,i\xi_n))<\infty.
$$
 Applications to systems of partial differential operators
      with constant coefficients are considered.
\end{abstract}

\section*{Introduction}

\subsection*{One-parameter semigroups in the convolution algebra of rapidly decreasing distributions}
Let $\Mmm$ be the set of $\mm$ matrices with complex entries, and $\calO^\prime_C(\Rn;\allowbreak\Mmm)$ the convolution algebra of $\Mmm$-valued distributions on $\Rn$ rapidly decreasing in the sense of L.~Schwartz. The Fourier transformation $\scrF$ is an isomorphism of $\calO^\prime_C(\Rn;\Mmm)$
onto the algebra $\calO_M(\Rn;\Mmm)$ of $\Mmm$-valued infinitely differentiable slowly increasing functions on~$\Rn$. We prove that $G\in\calO^\prime_C(\Rn;\Mmm)$ is the generating distribution of a one-parameter infinitely differentiable convolution semigroup $(S_t)_{t\ge 0}\subset \calO^\prime_C(\Rn;\Mmm)$
if and only if\vg
\begin{equation}
\max\{\hRe\lambda:\lambda\in\sigma((\scrF G)(\xi))\}= O(\log|\xi|)\ \quad\mbox{as }|\xi|\to\infty.\tag{i}
\end{equation}
In the above, $\sigma$ denotes the spectrum of a square matrix.

If $G=\calG(\partial_1,\ldots,\partial_n)\otimes\delta$ where $\delta$ is the Dirac distribution on~$\Rn$, $\partial_1,\ldots,\partial_n$ denote the first order partial derivatives with respect to the coordinates of~$\Rn$, and
$\calG(\partial_1,\ldots,\partial_n)$ is an $\mm$ matrix whose entries are scalar partial differential operators (PDOs) with constant coefficients,
then $(\scrF G)(\xi)=\calG(i\xi)$ for every $\xi\in\Rn$, and condition (i) takes the form
$$
\max\{\hRe\lambda:\lambda\in\sigma(\calG(i\xi))\}=O(\log|\xi|)
\ \quad
\mbox{as }|\xi|\to\infty.\eqno{{\rm(i)}^\prime}
$$
Thanks to the fact that $\det(\lambda\symjed_{\mm}-\calG(\zeta_1,\ldots,\zeta_n))$ is a polynomial, L.~G{\aa}r\-ding was able to prove the conjecture of I.~G.~Petrovski\u\i\  that (i)$'$ is equivalent to the condition
\begin{equation}
\sup\{\hRe\lambda:\lambda\in\sigma(\calG(i\xi)),\,\xi\in\Rn\}<\infty.\tag{ii}
\end{equation}

\subsection*{Application to the Cauchy problem for partial differential equations with constant coefficients}
Suppose that $\calG(\partial_1,\ldots,\partial_n)$ satisfies (ii), and
$(S_t)_{t\ge 0}\subset\calO^\prime_C(\Rn;\Mmm)$ is the infinitely differentiable convolution semigroup with generating distribution
$\calG(\partial_1,\ldots,\partial_n)\otimes\delta$. Suppose moreover that
$$
%\item[(iii)] {\spaceskip.33em plus.22em minus.17em
\textbot{$E$ is a sequentially complete l.c.v.s. continuously imbedded in $\calS^\prime(\Rn;\Cm)$ such that $(S_t \,*)E\subset E$
for every $t\in[0,\infty\[$, and the mapping $[0,\infty\[\times E\ni(t,u)\mapsto S_t*u\in E$
is separately continuous.}\eqno{\rm(iii)}
$$
Then $((S_t\,*)|_E)_{t\ge 0}\in L(E;E)$ is a one-parameter operator semigroup
of class $(C_{0})$ whose infinitesimal generator $\calG_E$ satisfies the equalities
\begin{align*}
D(\calG_E)&=\{u\in E:\calG(\partial_1,\ldots,\partial_n)u\in E\},\\
\calG_Eu&=\calG(\partial_1,\ldots,\partial_n)u\ \quad
\mbox{if }
u\in D(\calG_E).
\end{align*}
We prove that if (iii) holds, then for every $k=1,2,\ldots$ the Cauchy problem
\begin{equation}
\frac{d}{dt}u(t)=\calG(\partial_1,\ldots,\partial_n)u(t)\ \quad 
\mbox{for }
t\in[0,\infty\[,\ \quad u(0)=u_0, \tag{iv}
\end{equation}
with given $u_0\in D(\calG^k_E)$ has a solution $u(\cdot)\in C^k([0,\infty\[;E)$ which is unique in the class
$C^1([0,\infty\[;\calS^\prime(\Rn;\Cm))$. This solution is given by the formula
\begin{equation}
u(t)=S_t*u_0\ \quad\mbox{for }t\in[0,\infty\[.\tag{v}
\end{equation}
Examples of spaces $E$ satisfying (iii) are given in Sec.~8.

\subsection*{Hyperbolic partial differential systems with constant coefficients}
The matricial partial differential operator $\symjed_{\mm}\otimes\partial_t-\calG(\partial_1,\ldots,\partial_n)$ on $\symR^{1+n}=\{(t,x_1,\ldots,x_n)\}$ is called
{\it hyperbolic with respect to the coordinate} $t$ if (ii) holds and the hyperplane $t=0$ is non-characteristic for the operator. This last holds if and only if
$$
\textbot{the degree of the polynomial of $1+n$ variables
$$
P(\lambda,\zeta_1,\ldots,\zeta_n)=\det(\lambda\symjed_{\mm}-\calG(\zeta_1,\ldots,\zeta_n))
$$
is equal to $m$.}\eqno{\rm(vi)}
$$
Suppose that (ii) is satisfied and $(S_t)_{t\ge0}\subset\calO'_{C}(\Rn;\Mmm)$ is the infinitely differentiable convolution semigroup whose generating distribution is\break $\calG(\partial_1,\ldots,\partial_n)\otimes\delta$. Then the question arises about properties of
$(S_t)_{t\ge0}$
 corresponding to (vi). We prove that
\begin{aitem} 
\item[\rm(a)] if (vi) holds, then $(S_t)_{t\ge0}$ extends to a one-parameter convolution group $(S_t)_{t\in\symR}$ such that $\supp S_t$ is bounded for every  $t\in\symR$,
and 
\item[(b)] if (vi) does not hold, then $\supp S_t$ is unbounded for every $t\in\]0,\infty\[$.
\end{aitem}
	
\section{The setting and results}\label{sec1}
\subsection{Notation}\label{subsec1.1}
Throughout the present paper the symbols $\partial_1,\ldots,\partial_n$ denote partial derivatives of the first order (not multiplied by any constant) of a function or distribution on~$\Rn$. For partial derivatives of higher order we use the abbreviation
$\partial^\alpha=\partial^{\alpha_1}_1\ldots\partial^{\alpha_n}_n$ where
$\alpha=(\alpha_1,\ldots,\alpha_n)\in\symN_0^n$ is a {\it multiindex} whose {\it length} is defined as $|\alpha|=\alpha_1+\cdots+\alpha_n$. $\calS(\Rn)$ and $\calS^\prime(\Rn)$  denote the space of {\it infinitely differentiable rapidly decreasing complex functions} on~$\Rn$ and the space of {\it slowly increasing distributions} on $\Rn$. The Fourier transformation $\scrF$ is defined by the formulas
\begin{align}
(\scrF\varphi)(\xi_1,\ldots,\xi_n)&=\hat\varphi(\xi_1,\ldots,\xi_n)\notag\\
&=\mathop{\int\cdots\int}_{\Rn} e^{-i\sum^n_{k=1} x_k\xi_k}\varphi(x_1,\ldots, x_n)\,dx_1\ldots dx_n\label{eq1.1}
\end{align}
whenever $\varphi\in\calS(\Rn)$ and $(\xi_1,\ldots,\xi_n)\in\Rn$, and
\begin{equation}
\label{eq1.2}
\<\scrF T,\varphi\>=\<T,\scrF\varphi\>
\end{equation}
whenever $T\in\calS^\prime(\Rn)$, $\varphi\in\calS(\Rn)$ and $\scrF\varphi$ is determined by \eqref{eq1.1}. The compatibility of \eqref{eq1.2} with \eqref{eq1.1} follows from the Parseval equality for a pair of elements of $\calS(\Rn)$.

\subsection{The function algebra $\calO_M(\Rn)$ and the convolution algebra of distributions $\calO_C^\prime(\Rn)$}\label{subsec1.2}
 Let $\calO_M(\Rn)$ be the space of {\it infinitely differentiable slowly increasing} complex functions on~$\Rn$. Recall that $\phi\in\calO_M(\Rn)$ if and only if for every $\alpha\in\symN_0^n$ there is $m_\alpha\in\symN_0$
such that
$$ 
\sup_{\xi\in\Rn}(1+|\xi|)^{-m_\alpha}|\partial^\alpha\phi(\xi)|<\infty.
$$
Obviously $\calO_M(\Rn)$ is a function algebra. Furthermore
\begin{equation}
\label{eq1.3}
\calO_M(\Rn)=\{\phi\in C^\infty(\Rn):\phi\cdot\varphi\in \calS(\Rn)\ 
\mbox{for every } \varphi\in\calS(\Rn)\}.
\end{equation}

For every $k\in\symN_0$ denote by $\B_k(\Rn)$ the space of continuous complex functions
$f$ on $\Rn$ such that $f(x)=O(|x|^{-k})$ as $|x|\to\infty$. 
A distribution $T\in\calD^\prime(\Rn)$ is called {\it rapidly decreasing} if for every $k\in\symN_0$ there is $m_k\in\symN_0$ such that
$T=\sum_{|\alpha|\le m_k}\partial^\alpha f_{k,\alpha}$ where $f_{k,\alpha}\in\B_k(\Rn)$
for every $\alpha\in\symN_0^n$ with $|\alpha|\le m_k$. The space of rapidly decreasing distributions on~$\Rn$, denoted by $\calO'_C(\Rn)$, is a convolution algebra~$^{*)}$\footnote{$^{*)}$ See \cite[Sec. VII.5, pp.~246--248]{S3}, \cite[pp.~131--134]{K-R}.}. 
One has
\begin{equation}
\label{eq1.4}
\calO_C^\prime(\Rn)=\{T\in \calS^\prime(\Rn):T*\varphi\in\calS(\Rn)\ \mbox{for every }
\varphi\in\calS(\Rn)\}.
\end{equation}

{\spaceskip.33em plus.22em minus.17em Since $\calO_M(\Rn)$ and $\calO_C^\prime(\Rn)$ are subsets of $\calS^\prime(\Rn)$, 
$\scrF\calO_M(\Rn)$ and $\scrF\calO_C^\prime(\Rn)$
make sense. Furthermore,}
\begin{equation}
\label{eq1.5}
\scrF\calO_C^\prime(\Rn)=\calO_M(\Rn)
\end{equation}	
and
\begin{equation}
\label{eq1.6}
\scrF(U*V)=(\scrF U)\cdot(\scrF V)
\end{equation}	
for every $U,V\in\calO_C^\prime (\Rn)$~$^{*)}$\footnote{$^{*)}$ See \cite[Sec. VII.8,
Theorem XV, p.~268]{S3}, \cite[Theorem 8.23, p.~156]{K-R}.}.
The equality \eqref{eq1.6} means that the Fourier transformation is an (algebraic) isomorphism of the convolution algebra of distributions $\calO_C^\prime(\Rn)$ onto the function algebra $\calO_M(\Rn)$.

By the closed graph theorem, it follows from \eqref{eq1.3} and \eqref{eq1.4} that the operators $\phi\,\cdot$ for $\phi\in\calO_M(\Rn)$ and $T\,*$ for $T\in\calO_C^\prime(\Rn)$ belong to the space $L(\calS(\Rn);\calS(\Rn))$ of continuous linear operators from $\calS(\Rn)$ into $\calS(\Rn)$. Let $L_b(\calS(\Rn);\calS(\Rn))$ denote the space
$L(\calS(\Rn);\calS(\Rn))$ equipped with the compact-open topology. The sets of operators $\calO_M(\Rn)\,\cdot$ and $\calO_C^\prime(\Rn)\,*$ are closed subspaces of
$L_b(\calS(\Rn);\calS(\Rn))$, and we 
treat them as equipped with the induced topology. 
The Fourier transformation is a continuous isomorphism of
$\calO_C^\prime(\Rn)$ onto $\calO_M(\Rn)$. Furthermore, the bilinear maps
$\calO_M(\Rn)\times\calO_M(\Rn)\ni(\phi,\psi)\mapsto\phi\cdot\psi\in\calO_M(\Rn)$
 and $\calO_C^\prime(\Rn)\times\calO_C^\prime(\Rn)\ni(S,T)\mapsto S*T\in\calO_C^\prime(\Rn)$ are hypocontinuous. We shall prove the 
latter fact;
the proof of the former is the same. Since
$\calS(\Rn)$ is a barrelled space, the boundedness of a subset of
$L_b(\calS(\Rn);\calS(\Rn))$
 is equivalent to its equicontinuity. This implies that composition in $L_b(\calS(\Rn);\calS(\Rn))$ is hypocontinuous. Since for $U,V\in\calO_C^\prime(\Rn)$
one has $(S\,*)|_{\calS(\Rn)},(V\,*)|_{\calS(\Rn)}\in L(\calS(\Rn);\calS(\Rn))$ and
$((U*V)*)|_{\calS(\Rn)}=(U*)|_{\calS(\Rn)}\circ(V*)|_{\calS(\Rn)}$, it follows that convolution in $\calO_C^\prime(\Rn)$ is hypocontinuous.

\subsection{The function algebra  $\calO_M(\Rn;\Mmm)$ and the convolution algebra of distributions $\calO_C^\prime(\Rn;\Mmm)$}\label{subsec1.3} 
 Let $m,n\in\symN$, and let $\Mmm$ be the set of $\mm$ matrices with complex entries. Denote by $\calO_M(\Rn;\Mmm)$ the space of functions of the form
$\phi:\Rn\ni\xi\mapsto(\phi_{j,k}(\xi))^m_{j,k=1}\in\Mmm$ such that
$\phi_{j,k}\in\calO_M(\Rn)$ for all $j,k$.
This spaces carries the topology of $\calO_M(\Rn)^{m^2}$ where each factor is equipped with the topology induced by $L_b(\calS(\Rn);\calS(\Rn))$. Multiplication in
$\calO_M(\Rn;\Mmm)$ is defined by the rule
$$
(\phi\cdot\psi)(\xi)=\Big(\sum_{j=1}^m\phi_{i,j}(\xi)\psi_{j,k}(\xi)\Big)^m_{i,k=1}.
$$
$\calO_M(\Rn;\Mmm)$ is a locally convex algebra with hypocontinuous multiplication.

Denote by $\calO_C^\prime(\Rn;\Mmm)$ the space of $\mm$ matrices
$T=(T_{j,k})_{j,k=1}^m$ such that $T_{j,k}\in\calO_C^\prime(\Rn)$ for all
$j,k$. The convolution in $\calO_M(\Rn;\Mmm)$
is defined by the rule
$$
S*T=\Big(\sum_{j=1}^m S_{i,j}*T_{j,k}\Big)^m_{i,k=1}.
$$
The space $\calO_C^\prime(\Rn;\Mmm)$
 carries the topology of $\calO_C^\prime(\Rn)^{m^2}$ where each factor is equipped with the topology induced by $L_b(\calS(\Rn);\calS(\Rn))$. The l.c.v.s. $\calO_C^\prime(\Rn;\Mmm)$ is a locally convex associative convolution algebra of $\Mmm$-valued distributions.
Convolution in $\calO_C^\prime(\Rn;\Mmm)$ is hypocontinuous.

The analogues of \eqref{eq1.5} and \eqref{eq1.6} are valid for
$\calO_M(\Rn;\Mmm)$ and\break $\calO_C^\prime(\Rn;\Mmm)$.

\subsection{{\spaceskip.33em plus.22em minus.17em Infinitely differentiable one-parameter convolution semigroups in $\calO_C^\prime(\Rn;\Mmm)$}}\label{subsec1.4}
 By a one-parameter {\it infinitely differentiable convolution semigroup in}\break $\calO_C^\prime(\Rn;\Mmm)$, briefly 
{\it i.d.c.s.}, we mean a mapping
\begin{equation}
\label{eq1.7}
[0,\infty\[\ni t\mapsto S_t\in\calO_C^\prime(\Rn;\Mmm)
\end{equation}
such that
%\bnumer
%\item
\begin{equation}\label{eq1.8}
\textbot{$S_{s+t}=S_s*S_t$ for every $s,t\in[0,\infty\[$,}
\end{equation}
%\item
\begin{equation}\label{eq1.9}
\textbot{$S_0=\symjed_{\mm}\otimes\delta$ where $\symjed_{\mm}$  is the unit $\mm$ matrix and
$\delta$  is the Dirac distribution on~$\Rn$,}
\end{equation}
%\enumer\vskip-\lastskip\vskip0pt
%\dbnumer
%\item
\begin{equation}\label{eq1.10} 
\textbot{\hspace{3pt}the mapping \eqref{eq1.7} is infinitely differentiable.}
\end{equation}
%\enumer
In \eqref{eq1.10} it is understood that the derivatives at zero are right derivatives, and that the topology in $\calO_C^\prime(\Rn;\Mmm)$ is that defined in Sec.~\ref{subsec1.3}.

The {\it generating distribution} of the i.d.c.s. $(S_t)_{t\ge 0}
\subset\calO_C^\prime(\Rn;\Mmm)$ is defined as
$$
G:=\frac d{dt}\bigg|_{t=0}S_t\in\calO_C^\prime(\Rn;\Mmm).
$$
It follows that
$$
\frac d{dt}S_t=G*S_t=S_t*G\ \quad\mbox{for every }t\in[0,\infty\[.
$${\spaceskip.33em plus.22em minus.17em 
Furthermore, {\it any i.d.c.s. in $\calO_C^\prime(\Rn;\Mmm)$ is uniquely determined by its generating distribution}. Indeed, suppose that $G\in\calO_C^\prime(\Rn;\allowbreak\Mmm)$ is the generating distribution of two i.d.c.s. $(S_t)_{t\ge 0},(T_{t})_{t\ge 0}\!\subset\!\calO_C^\prime(\Rn;\allowbreak\Mmm)$. Fix any
$t\in\]0,\infty\[$. Then
$(S_\tau\,*)|_{\calS(\Rn;\Cm)}, (T_{t-\tau}\,*)|_{\calS(\Rn;\Cm)}\in L(\calS(\Rn;\allowbreak\Cm);
\calS(\Rn;\Cm))$ and}
$$
((S_\tau*T_{t-\tau})*)|_{\calS(\Rn;\Cm)}=(S_\tau\,*)|_{\calS(\Rn;\Cm)}\circ
(T_{t-\tau}*)|_{\calS(\Rn;\Cm)} \quad\mbox{for every }\tau\!\in\![0,t].
$$
Since $\calS(\Rn;\Cm)$ is a Montel (and hence barrelled) space, one infers from the Banach--Steinhaus theorem that the function
$[0,t]\ni\tau\mapsto S_\tau*T_{t-\tau}\in\calO_C^\prime(\Rn;\Mmm)$ 
is continuously differentiable and
$$
\frac d{d\tau}(S_\tau*T_{t-\tau})=\bigg(\frac d{d\tau}S_\tau\bigg)*T_{t-\tau}
+S_\tau*\bigg(\frac d{d\tau}T_{t-\tau}\bigg). 
$$
Consequently,
$$ 
\frac d{d\tau}(S_\tau*T_{t-\tau})=(S_\tau* G)*T_{t-\tau}-S_\tau*(G*T_{t-\tau})
=0, 
$$
by associativity of the convolution in
$\calO_C^\prime(\Rn;\Mmm)$, so that $S_\tau*T_{t-\tau}$  is independent of
$\tau$ for $\tau\in[0,t]$, and 
$S_t=(S_\tau*T_{t-\tau})|_{\tau=t}=(S_\tau*T_{t-\tau})|_{\tau=0}=T_{t}$.

The Cauchy problem for a PDO with constant coefficients can be reduced by Fourier transformation with respect to the spatial coordinates to the Cauchy problem with a parameter for an ODO. In the framework of the spaces
$\calO_C^\prime(\Rn;\Mmm)$
and $\calO_M(\Rn;\Mmm)$ this method consists in making use of the following

\begin{xlem} Suppose that $G\in\calO_C^\prime(\Rn;\Mmm)$ and let
$A=\scrF G$, so that $A\in\calO_M(\Rn;\Mmm)$. Then the following two conditions are equivalent:
\begin{aitem}
\item[\rm(a)] $G$ is the generating distribution of the i.d.c.s. $(S_t)_{t\ge0}
\subset\calO_C^\prime(\Rn;\Mmm)$,
\item[\rm(b)] $\exp (tA(\cdot))\in\calO_M(\Rn;\Mmm)$ for every $t\in[0,\infty\[$ and the mapping $[0,\infty\[\ni t\mapsto \exp (tA(\cdot))\in\calO_M(\Rn;\Mmm)$
is infinitely differentiable.
\end{aitem}
Furthermore, if $A=\scrF G$  and {\rm(a)}, {\rm(b)} are satisfied, then
$\exp (tA(\cdot))=\scrF S_t$ and
$$
(S_t\,*)|_{\calS^\prime(\Rn;\Cm)} =\scrF^{-1}\circ[(\exp tA(\cdot))\,\cdot]
\circ\scrF|_{\calS^\prime(\Rn;\Cm)}\ \quad
\mbox{for every }t\in[0,\infty\[.
$$
\end{xlem}

Basing on the above lemma we shall prove four theorems. For this purpose we shall use some intricate facts concerning $\calO_C^\prime$ and $\calO_M$, which for the most part are only mentioned in \cite{S3}, and are presented in detail in \cite{K3}. For any 
$B\in\Mmm$
denote by $\sigma(B)$ the spectrum of the matrix~$B$.

\begin{thm}\label{thm1} A distribution $G\in\calO_C^\prime(\Rn;\Mmm)$ is the generating distribution of an i.d.c.s.
$(S_t)_{t\ge0}\subset\calO_C^\prime(\Rn;\Mmm)$
if and only if
\begin{equation}
\label{eq1.11}
\max\{\hRe\lambda:\lambda\in\sigma((\scrF G)(\xi))\}=O(\log|\xi|)\ \quad\mbox{as }|\xi|\to\infty,\,
\xi\in\Rn.
\end{equation}
\end{thm}

The quantity
\begin{equation}
\label{eq1.12}
s(G):=\sup\{\hRe\lambda:\mbox{there is }\xi\in\Rn\mbox{ such that }
\lambda\in\sigma((\scrF G)(\xi))\},
\end{equation}
finite or equal to $+\infty$, will be called the {\it spectral bound} of~$G$. For any i.d.c.s. $(S_t)_{t\ge0}\subset\calO_C^\prime(\Rn;\Mmm)$ let
%\dbnumer
%\item
\begin{equation}\label{eq1.13}
\textbot{$\omega((S_t)_{t\ge0}):=\inf\{\omega\in\Rn:\mbox{the}$ one-parameter semigroup of operators 
$$
((e^{-\omega t}S_t\,*)|_{\calS(\Rn;\Cm)})_{t\ge 0}\subset L(\calS(\Rn;\Cm);
\calS(\Rn;\Cm))
$$
is equicontinuous\}
}
\end{equation}
%\enumer
where it is assumed that $\inf\emptyset=+\infty$. We call
$\omega((S_t)_{t\ge0})$ the {\it growth bound} of the i.d.c.s. $(S_t)_{t\ge0}$~$^{*)}$\footnote{$^{*)}$ In \eqref{eq1.13} the {\it growth bound with respect to} $\calS(\Rn;\Cm)$ is defined. The growth bounds with respect to some other spaces invariant for the semigroup $(S_t\,*)_{t\ge0}$ are also
equal to the spectral bound of the generating distribution.
See \cite{B} and \cite[Theorem 1]{K2}. For one-parameter semigroups of operators in a Banach space the relations between
the growth bound of the semigroup and the spectral bound of its generator are discussed in great detail in \cite[Sec.~IV.2]{E-N}.}.

\begin{thm}\label{thm2} For every i.d.c.s. $(S_t)_{t\ge0}\subset\calO_C^\prime(\Rn;\Mmm)$
its growth bound is equal to the spectral bound of its generating distribution.
\end{thm}

Let $\calG(\partial_1,\ldots,\partial_n)$ be an $\mm$ matrix whose entries are PDOs on $\Rn$ with constant complex coefficients. Let $\delta$ be the Dirac distribution on~$\Rn$. Then $\calG(\partial_1,\ldots,\partial_n)\otimes\delta\in\calO_C^\prime(\Rn;\Mmm)$
and
$$
[\scrF(\calG(\partial_1,\ldots,\partial_n)\otimes\delta)](\xi)=\calG(i\xi_1,\ldots,i \xi_n)
$$
for every $\xi=(\xi_1,\ldots, \xi_n)\in\Rn$. The quantity
$$
s_0(\calG):=\sup\{\hRe\lambda:\lambda \in\sigma(\calG(i\xi_1,\ldots,i \xi_n))
,\,(\xi_1,\ldots,\xi_n)\in\Rn\}
$$
is equal to the spectral bound of the distribution $\calG(\partial_1,\ldots,\partial_n)\otimes\delta$. It was conjectured by I.~G.~Petrovski\u\i\  
\cite[footnote on p.~24]{P} and proved by L.~G{\aa}rding \cite[Lemma on p.~11]{G} that $s_0(\calG)<\infty$ if and only if
$G=\calG(\partial_1,\ldots,\partial_n)\otimes\delta$ satisfies \eqref{eq1.11}~$^{**)}$\footnote{$^{**)}$ If $G=\calG(\partial_1,\ldots,\partial_n)\otimes\delta$,
then \eqref{eq1.11} takes the form $\sup\{\hRe\lambda:\lambda\in\sigma(\calG(i\xi))\}
=O(\log|\xi|)$ as $|\xi|\to\infty$, and in this form \eqref{eq1.11} occurs in \cite[Sec.~I.5]{P}. However, usually the ~~``Petrovski\u\i\ condition'' means the assumption that $s_0(\calG)<\infty$.}. 
Therefore Theorems 1.1 and 1.2 imply

\begin{thm}\label{thm3} Let $\calG(\partial_1,\ldots,\partial_n)$ be an $\mm$ matrix whose entries are PDOs on $\Rn$ with constant complex coefficients. Then the following two conditions are equivalent:
%\dbnumer
%\item
\begin{gather}\label{eq1.14}
\textbot{$s_0(\calG)<\infty$,}\\
%
%\item
%\begin{equation}
\label{eq1.15} 
\textbot{$\calG(\partial_1,\ldots,\partial_n)\otimes\delta$ is the generating distribution\vadjust{\vskip2pt} of an i.d.c.s. $(S_t)_{t\ge0}\subset\calO_C^\prime(\Rn;\Mmm)$.}
\end{gather}
%\enumer
Furthermore, if these equivalent conditions are fulfilled, then there is exactly one i.d.c.s.
$(S_t)_{t\ge0}\subset\calO_C^\prime(\Rn;\Mmm)$ satisfying \eqref{eq1.15}, and the growth bound of this i.d.c.s.
 is equal to $s_0(\calG)$.
\end{thm}

\begin{exa}\label{exa1} Let $m=1$, $\calG(\partial_1,\ldots,\partial_n)=i(\partial_1^2+\cdots+\partial_n^2)$. Then $s_0(\calG)=s_0(-\calG)=0$, so that
$\calG(\partial_1,\ldots,\partial_n)\delta$ and $-\calG(\partial_1,\ldots,\partial_n)\delta$ are generating distributions of i.d.c.s. imbedded in $\calO_C^\prime(\Rn)$. Consequently, $i(\partial_1^2+\cdots+\partial_n^2)\delta$
is the generating distribution of an {\it infinitely differentiable one-parameter convolution group} $(S_t)_{t\in\Rn}\subset\calO_C^\prime(\Rn)$. This group of distributions satisfies the Schr\"odinger partial differential equation
$$
\partial_t S_t=i(\partial_1^2+\cdots+\partial_n^2)S_t,
$$
one has $S_0=\delta$, and for every $t\in\Rn\setminus\{0\}$ the distribution
$S_t$ is equal to the bounded function belonging to $\calO_M(\Rn)$ such that
$$
S_t(x)=(4\pi it)^{-n/2}\exp\bigg(\frac{i|x|^2}{4t}\bigg)\ \quad
\mbox{whenever }x\in\Rn.
$$
The factor $(4\pi it)^{-n/2}$ is defined as $\big(\frac 1{\sqrt{4\pi it}}\big)^n$ where $\arg\sqrt{4\pi it}=(\pi/4)\sgn t$. See \cite[p.~54]{Go}, \cite[p.~107]{R},
\cite[p.~48]{S1}. The direct proof that $S_t\subset\calO_C^\prime(\Rn)$ is on p.~245 of \cite{S2}. Another proof is by Fourier transformation: one has
$(\scrF S_t)(\xi)=e^{-it|\xi|^2}$, so that $\scrF S_t\in \calO_M(\Rn)$, and hence
$S_t\subset \calO_C^\prime(\Rn)$.
\end{exa}

\begin{exa}\label{exa2} Following J.~Rauch \cite[Sec.~3.10]{R} we look for solutions of class $C^\infty([0,\infty\[;\calS(\Rn))$ of the Cauchy problem
\begin{multline}\label{eq1.16}
\kern\parindent\sum_{k=0}^m
Q_k(\partial_1,\ldots,\partial_n)\partial_t^k u(t,x_1,\ldots,x_n)
=0\\	 
\hfill\mbox{for }(t,x_1,\ldots,x_n)\in[0,\infty\[\times\Rn,\kern\parindent\phantom{(1.16)}\\
\kern\parindent\partial_t^ku(0,x_1,\ldots,x_n)=u_k(x_1,\ldots,x_n)\hfill\\
\hfill\mbox{ for }k=0,\ldots,m-1
\mbox{ and }(x_1,\ldots,x_n)\in\Rn,
\end{multline}
where $Q_k(\partial_1,\ldots,\partial_n)$, $k=0,\ldots,m$, are linear partial differential operators with constant coefficients, and $u_k\in\calS(\Rn)$, $k=0,\ldots,m-1$, are given. As in \cite[Sec.~3.10]{R}, we assume that the polynomial
$$
P(\lambda,\zeta_1,\ldots,\zeta_n)=Q_m(\zeta_1,\ldots,\zeta_n)\lambda^m+\cdots+
Q_1(\zeta_1,\ldots,\zeta_n)\lambda+ Q_0(\zeta_1,\ldots,\zeta_n)
$$
has two properties:
%\dbnumer
%\item
\begin{gather}\label{eq1.17} 
\textbot{$\sup\{\hRe\lambda:\lambda\in\symC$, there is $\xi\in\Rn$ such that
$P(\lambda,i\xi)=0\}=s_0<\infty$,}\\
%\item
\label{eq1.18}  
\textbot{$Q_m(i\xi)\ne0$ whenever $\xi\in\Rn$.}
\end{gather}
%\enumer
For every $\xi\in\Rn$ denote by $\hat G(\xi)$ the matrix
$$
\left[\begin{matrix}
0&1&&&\\
&0\\
&&&1&\\
&&&0&1\\
-\frac{Q_0(i\xi)}{Q_m(i\xi)}&-\frac{Q_1(i\xi)}{Q_m(i\xi)}&\cdots&-\frac{Q_{m-2}(i\xi)}{Q_m(i\xi)}&-\frac{Q_{m-1}(i\xi)}{Q_m(i\xi)}
\end{matrix}\right].
$$
By \cite[Example A.2.7]{H} there is $m_0\in\symN$ such that
$$
\sup\{(1+|\xi|)^{-m_0}|Q_m(i\xi)^{-1}|:\xi\in\Rn\}<\infty.
$$
Since $\partial^\alpha(Q_m(i\xi)^{-1})=Q_m(i\xi)^{-1-|\alpha|}R(\xi)$ for every
$\alpha\in\symN_0^n$ where $R(\xi)$ is a polynomial, it follows that
$Q_m(i\,\cdot)^{-1}\in\calO_M(\Rn)$. Consequently,
\begin{equation}
\label{eq1.19}
\hat G\in\calO_M(\Rn;\Mmm).
\end{equation}
Furthermore,
$$
\det(\lambda\symjed_{\mm}-\hat G(\xi))=\lambda^m+\frac{Q_{m-1}(i\xi)}{Q_m(i\xi)}\lambda^{m-1}
+\cdots+\frac{Q_1(i\xi)}{Q_m(i\xi)}\lambda+\frac{Q_0(i\xi)}{Q_m(i\xi)},
$$
whence
$$
\sigma(\hat G(\xi))=\{\lambda\in\symC:P(\lambda,i\xi)=0\},
$$
and so \eqref{eq1.17} implies that
\begin{equation}
\label{eq1.20}
\sup\{\hRe\lambda:\lambda\in\sigma(\hat G(\xi)),\,\xi\in\Rn\}=s_0.
\end{equation}
{\spaceskip.33em plus.22em minus.17em From \eqref{eq1.19} it follows that there is a unique distribution
$G\in\calO_C^\prime(\Rn;\Mmm)$ such that
$\scrF G=\hat G$. By Theorems \ref{thm1} and \ref{thm2}, \eqref{eq1.20} implies that
$G$ is the generating distribution of an i.d.c.s. $(S_t)_{t\ge 0}\subset\calO_C^\prime(\Rn;\Mmm)$ such that $\omega((S_t)_{t\ge 0})\allowbreak=s_0$ and}
$$
\left[\begin{matrix}
1\\
&\ddots\\
&&1\\
&&&Q_m
\end{matrix}\right]\partial_t S_t=
\left[\begin{matrix}
0&1\\
&0\\
&&&1\\
&&&0&1\\
-Q_0&-Q_1&\cdots&-Q_{m-2}&-Q_{m-1}
\end{matrix}\right] S_t
$$
for every $t\in[0,\infty\[$ where $Q_k=Q_k(\partial_1,\ldots,\partial_n)$ for $k=0,\ldots,m$. By arguments similar to that presented in Sec.~8, the above implies that, under the assumptions \eqref{eq1.17} and \eqref{eq1.18}, for every
$u_0,\ldots, u_{m-1}\in\calS(\Rn)$ and $u\in C^\infty([0,\infty\[;\calS(\Rn))$
the following two conditions are equivalent:
\begin{aitem}
\item[(a)] $u$ is a solution of the Cauchy problem \eqref{eq1.16},\vskip4pt
\item[(b)] $\left[\begin{matrix}
u(t,\cdot)\\
\partial_t u(t,\cdot)\\
\vdots\\
\partial_t^{m-1} u(t,\cdot)
\end{matrix}\right]
=
S_t*\left[\begin{matrix}
u_0\\
u_1\\
\vdots\\
u_{m-1}
\end{matrix}\right]
\mbox{ for }t\in[0,\infty\[$.
\end{aitem}

If only the condition \eqref{eq1.17} is satisfied and \eqref{eq1.18} may fail, then the i.d.c.s.'s in $\calO_C^\prime(\Rn;\Mmm)$ seem not to be useful, but
 $\calO_C^\prime(\symR^{1+n})$ can be used to express the properties of the fundamental solution for the operator $Q_m(\partial_1,\ldots,\partial_n)\partial_t^m+\cdots+
Q_1(\partial_1,\ldots,\partial_n)\partial_t+Q_0(\partial_1,\ldots,\partial_n)$
with support contained in $H_+=\{(t,x_1,\ldots,x_n)\in\symR^{1+n}:t\ge 0\}$. 
See the article of the present author in arXiv:1105.0877.
\end{exa}

\noindent{\bf Comments.} I.~G.~Petrovski\u\i\  \cite{P} was the first to notice the significance of smooth slowly increasing functions in the theory of evolutionary PDEs. 
The theory of distributions did not yet exist in 1938 when [P] was published, and only in 1950 did L.~Schwartz explain in \cite{S1} how the results of Petrovski\u\i\  may be elucidated by placing them in the framework of $\calO_C^\prime$. However in \cite{S1} the spectral properties of $[\scrF(\calG(\partial_1,\ldots,\partial_n)\otimes\delta)](\xi)
=\calG(i\xi,\ldots,i\xi_n)$ were not discussed.

If $G\!\in\!\calO_C^\prime(\Rn;\Mmm)$ and $\calN_G\!=\!\{(\lambda,\xi)\in\symC\times\Rn:
{\det(\lambda\symjed_{\mm}\!-\!(\scrF G)(\xi))}\allowbreak=0\}$, then \eqref{eq1.11} may be expressed in an equivalent form: there is $C\in\]0,\infty\[$ such that
$$
\mbox{if}\quad (\lambda,\xi)\in\calN_G,\quad\mbox{then}\quad
\hRe\lambda\le C(1+\log(1+|\xi|)).\eqno{(1.11)'}
$$
As mentioned earlier, just this logarithmic condition was used in \cite{P}. 
In connection with convolution equations similar logarithmic estimates (in $\symC^{1+n}$ instead of 
$\symC\times\Rn$) were used by L.~Ehrenpreis in \cite{E1} and in \cite[Sec.~VIII.3]{E2}. Logarithmic estimates related to convolution equations also occur in elaborate theorems of L.~H\"ormander \cite[Secs.~16.6 and 16.7]{H}. The role of conditions \eqref{eq1.14} and \eqref{eq1.11} in the theory of evolutionary PDOs with constant coefficients is discussed in \cite[Sec.~3.10]{R}. 

From the above-mentioned Petrovski\u\i\  conjecture proved by G{\aa}rding, and from Theorem~\ref{thm3}, it follows that whenever the generating distribution
$G\in\calO_C^\prime(\Rn;\Mmm)$
of an i.d.c.s. $(S_t)_{t\ge0}\subset\calO_C^\prime(\Rn;\Mmm)$  has the form
$G=\calG(\partial_1,\ldots,\partial_n)\otimes\delta$, then
$$
\sup\{\hRe\lambda:\lambda\in\sigma(\calG(i\xi_1,\ldots,i\xi_n)),\,{(\xi_1,\ldots,\xi_n)\in\Rn}\}
=s_0<\infty,
$$
and whenever $\varepsilon>0$, then the semigroup of operators
$$
((e^{-(s_0+\varepsilon)t}S_t\,*)|_{\calS(\Rn;\Cm)})_{t\ge0}\subset L(\calS(\Rn;\Cm);
\calS(\Rn;\Cm))
$$
is equicontinuous. As noticed by L.~Schwartz \cite{S2}, the theory of equicontinuous one-parameter semigroups of operators in an l.c.v.s. imitates the theory of one-parameter semigroups of operators in a Banach space. A detailed presentation of the theory of equicontinuous one-parameter semigroups of operators in a sequentially complete l.c.v.s. is contained in Chapter IX of the monograph of K.~Yosida \cite{Y}.

\subsection{Relation to hyperbolic systems of PDOs}\label{subsec1.5}
Let $\calE^\prime(\Rn)$ be the space of distributions on $\Rn$ with compact support,
equipped
with the topology of uniform convergence on bounded subsets of $C^\infty(\Rn)$.
L.~Ehrenpreis \cite[Sec.V.5]{E2} proved that $\calE^\prime(\Rn)=\{T\in\calD^\prime(\Rn):
T\,*\in L(\calD(\Rn);\calD(\Rn))\}$ and the topology induced in $\calE^\prime(\Rn)$  by
$L_b(\calD(\Rn);\allowbreak\calD(\Rn))$ via the mapping $T\mapsto T\,*$ coincides with the original topology of $\calE^\prime(\Rn)$. This topology is stronger than the one iduced on
$\calE^\prime(\Rn)$ by $\calO_C^\prime(\Rn)$. See \cite[Sec.~III.7]{S3}, \cite[Sec.~V.5, Lemma 5.17]{E2}. Let  $\calE^\prime(\Rn;\Mmm)$ be the space of $\Mmm$-valued distributions on $\Rn$ with compact support, i.e. the space of $\mm$ matrices whose entries belong to
$\calE^\prime(\Rn)$. With the topology of $\calE^\prime(\Rn)^{m^2}$ and convolution defined as in  $\calO_C^\prime(\Rn;\Mmm)$, the space $\calE^\prime(\Rn;\Mmm)$ is a convolution algebra with continuous convolution. See \cite[Sec.~VII.3, Theorem IV]{S3}. 

As in Theorem 1.3, let $\calG(\partial_1,\ldots,\partial_n)$ be an $\mm$ matrix whose entries are PDOs on $\Rn$ with constant complex coefficients. Put
\begin{equation}
\label{eq1.21}
P(\lambda,\zeta_1,\ldots,\zeta_n)=\det(\lambda\symjed_{\mm}-\calG(\zeta_1,\ldots,\zeta_n))
\end{equation}	
where $(\lambda,\zeta_1,\ldots,\zeta_n)\in\symC^{1+n}$.

\begin{thm}\label{thm4} Assume that $\calG(\partial_1,\ldots,\partial_n)$ satisfies condition \eqref{eq1.14}, and let $(S_t)_{t\ge0}\subset\calO_C^\prime(\Rn;\Mmm)$
be the i.d.c.s. whose generating distribution is
$\calG(\partial_1,\ldots,\partial_n)\otimes\delta$. Then the following three conditions are equivalent:
%\dbnumer
%\item
\begin{gather}\label{eq1.22} 
\textbot{there is $t_0\in\]0,\infty\[$ such that $S_{t_0}\in\calE'(\Rn;\Mmm)$,}\\[3pt]
\label{eq1.23} 
\textbot{the polynomial $P(\lambda,\zeta_1,\ldots,\zeta_n)$ defined by \eqref{eq1.21} has degree~$m$,}\\[3pt]	
\label{eq1.24}  %{\spaceskip.29em plus.22em minus.17em 
\textbot{$(S_t)_{t\ge 0}$ is an i.d.c.s. in the topological convolution\vadjust{\vskip2pt} algebra
$\calE^\prime(\Rn;\Mmm)$, and
may be uniquely extended\vadjust{\vskip2pt} to a one-parameter infinitely differentiable subgroup 	of $\calE^\prime(\Rn;\Mmm)$.}
\end{gather}
%\enumer	
\end{thm}

Let
\begin{equation}
\label{eq1.25}
N=\{(\lambda,\zeta_1,\ldots,\zeta_n)\in\symC^{1+n}:P(\lambda,\zeta_1,\ldots,\zeta_n)=0\}. 
\end{equation}
The matricial PDO
\begin{equation}
\label{eq1.26}
\symjed_{\mm}\otimes\partial_t-\calG(\partial_1,\ldots,\partial_n)
\end{equation}	
 on $\symR^{1+n}=\{(t,x_1,\ldots,x_n):t\in\symR,(x_1,\ldots,x_n)\in\Rn\}$
is said to be {\it hyperbolic in the sense of Ehrenpreis with respect to the coordinate} $t$ if there is $C\in\]0,\infty\[$ such that
\begin{equation}
\label{eq1.27}
\mbox{if } (\lambda,\zeta_1,\ldots,\zeta_n)\in N,\mbox{ then }
|\!\hRe\lambda|\le C(1+|\!\hRe\zeta_1|+\cdots+|\!\hRe\zeta_n|).
\end{equation}
Condition \eqref{eq1.27} is stronger than \eqref{eq1.14} which is equivalent to the existence of $C\in\]0,\infty\[$ 
such that
$$\displaylines{(1.14)^\prime\ \hfill
\mbox{if } (\lambda,\zeta_1,\ldots,\zeta_n)\in N\mbox{ and }\hRe\zeta_1=\cdots=\hRe\zeta_n=0,\mbox{ then }\hRe\lambda\le C.
\cr}$$
The matricial PDO \eqref{eq1.26} is said to be {\it hyperbolic in the sense of G{\aa}rding with respect to the coordinate} $t$ if the polynomial \eqref{eq1.21} satisfies \eqref{eq1.14}$^\prime$ and \eqref{eq1.23}. In the proof of Theorem \ref{thm4} it will be shown that for the matricial PDO \eqref{eq1.26} these two notions of hyperbolicity with respect to $t$ are equivalent.
Therefore
Theorem \ref{thm4} may be reformulated as follows: {\it if $\calG(\partial_1,\ldots,\partial_n)$
satisfies the Petrovski\u\i\ condition {\rm(1.14)}, then for the semigroup 
$(S_t)_{t\ge0}\!\subset\! \calO_C^\prime(\Rn;\Mmm)$ with generating distribution 
$\calG(\partial_1,\ldots,\partial_n)\otimes\delta$ the properties {\rm(1.22)} and {\rm(1.24)}
are equivalent, and they both hold if and only if the matricial PDO {\rm(1.26)} is hyperbolic with
respect to the variable}~$t$.

Suppose that \eqref{eq1.26} is hyperbolic with respect to~$t$. Let $P_m$ be the principal homogeneous part of the polynomial \eqref{eq1.21}, and let $\varGamma$ be the connected component of the set
$\{(\sigma,\xi_1,\ldots,\xi_n)\in\symR^{1+n}:P_m(\sigma,\xi_1,\ldots,\xi_n)\ne0\}$
which contains $(1,0,\ldots,0)$. By 
\cite[Lemma 8.7.3]{H}, $\varGamma$ is a convex cone. Let $\varGamma^0$ be the closed cone dual to~$\varGamma$. Using \cite[Theorem 12.5.1]{H} it may be proved that
\begin{equation} 
\label{eq1.28}
\varGamma^0=\{(t,x_1,\ldots,x_n)\in\symR^{1+n}:t\ge0,\, (x_1,\ldots,x_n)\in\convsupp S_t\}
\end{equation}	
where $(S_t)_{t\ge 0}$ is the i.d.c.s. occurring in Theorem~\ref{thm4}.
By \eqref{eq1.28}, the distribution $N\in\calD'(\Rn;\Mmm)$ such that
$\langle N,\varphi\rangle=\int_0^\infty\langle S_t,\varphi(t,\cdot)\rangle\,dt$
for every $\varphi\in\calD(\symR^{1+n})$ is a fundamental solution of (1.26) with support contained in~$\Gamma^0$. Theorem 4 resembles Theorems V and VI of \cite[Sec.~13]{S1}, and Theorems 12.5.1 and 12.5.2 of \cite{H}.

\section{A link between properties of  $\Mmm$-valued functions $\xi\mapsto A(\xi)$
and $(t,\xi)\mapsto\exp(tA(\xi))$}\label{sec2}

\begin{thms}[The Shilov inequality]\label{thms2.1}
Let $A\in\Mmm$. Then for every $t\in[0,\infty\[$ one has
\begin{equation}
\label{eq2.1}
\|\!\exp(t A)\|_{\Mmm}\le\rho(\exp(t A))\bigg(1+\sum_{k=1}^{m-1}
\frac{(2t)^k}{k!}\|A\|_{\Mmm}^k\bigg)
\end{equation}	
and
\begin{equation}
\label{eq2.2}
\rho(\exp(t A))=e^{t\max\hRe\sigma(A)}
\end{equation}	
where $\rho$ stands for the spectral radius, and $\sigma(A)$ denotes the spectrum of~$A$.
\end{thms}

The equality \eqref{eq2.2} follows from the spectral mapping theorem. The Shilov inequality \eqref{eq2.1} is an elaborate result of the theory of functions of matrices. See \cite{Sh}, \cite[Sec.~I.4]{Ge}, \cite[Sec.~II.6]{G-S}, 
\cite[Sec.~7.2]{F}. We say that $\Phi\subset C^\infty(\Rn;\Mmm)$
is a {\it set of uniformly slowly increasing functions} if for every 
$\alpha\in\symN_0^n$ there is $k_\alpha\in\symN_0$ such that
$\sup\{(1+|\xi|)^{-k_\alpha}\|(\partial/\partial\xi)^\alpha\phi(\xi)\|_{\Mmm}:\phi\in\Phi,\,
\xi\in\Rn\}<\infty$.

\begin{props}\label{props2.2} For any $A(\cdot)\in\calO_M(\Rn;\Mmm)$ the following three conditions are equivalent:
%\bnumer
\begin{gather}\label{eq2.3} 
\textbot{$\max\hRe\sigma(A(\xi))=O(\log|\xi|)$ as $|\xi|\to\infty$,}\\[3pt]
\label{eq2.4} \textbot{for every $T\in\]0,\infty\[$ there are $C\in\]0,\infty\[$ and
$k\in\symN$ such that
$$
\|\!\exp(tA(\xi))\|_{\Mmm}\le C(1+|\xi|)^k
$$
whenever 
$t\in[0,T]$ and $\xi\in\Rn$,}\\[3pt]
\label{eq2.5} 
\textbot{whenever $T\in\]0,\infty\[$, then $\{\exp(tA(\cdot)):
t\in[0,T]\}$ is a set of uniformly slowly increasing infinitely differentiable
$\Mmm$-valued functions on~$\Rn$.}
\end{gather}
\end{props}

\begin{props}\label{props2.3} For every $A(\cdot)\in\calO_{M}(\Rn;\Mmm)$ and $s_{0}\in\symR$ the following five conditions are equivalent:
\begin{gather}
\label{eq2.6} 
\textbot{$\sup\{\hRe\lambda:\lambda\in\sigma(A(\xi)),\,\xi\in\Rn\}\le s_0$;}\\[10pt]
\label{eq2.7} 
\textbot{there is $k\in\symN_0$ such that for every $\varepsilon>0$,
$$
\displaylines{\sup\{e^{-(s_0+\varepsilon)t}(1+|\xi|)^{-k}\|\!\exp(tA(\xi))\|_{\Mmm}:\hfill\cr
\hfill t\in[0,\infty\[,\,
\xi\in\Rn\}<\infty;}
$$}
\end{gather}\vspace{4pt}
$$
\textbot{for every $\varepsilon>0$ there is $k\in\symN$  such that
$$
\displaylines{
\sup\{e^{-(s_0+\varepsilon)t}(1+|\xi|)^{-k}\|\!\exp(tA(\xi))\|_{\Mmm}:\hfill\cr
\hfill t\in[0,\infty\[,\,
\xi\in\Rn\}<\infty;}
$$}\eqno{(2.7)^*}
$$\vspace{4pt}
\begin{equation}
\label{eq2.8} 
\textbot{for every $\alpha\in\symN_0^n$ there is $k_\alpha\in\symN_0$ such that
for every $\varepsilon>0$,
$$\displaylines{
\indent \sup\{e^{-(s_0+\varepsilon)t}(1+|\xi|)^{-k_{\alpha}}\|(\partial/\partial\xi)^{\alpha}\exp(tA(\xi))\|_{\Mmm} :\hfill\cr 
\hfill t\in\{0,\infty\[,\,
\xi\in\Rn\}<\infty;
}
$$}
\end{equation}\vspace{4pt}
$$
 \textbot{whenever $\varepsilon\in\]0,\infty\[$, then $\{e^{-(s_0+\varepsilon)t}
\exp(tA(\cdot)):t\in[0,\infty\[\}$ is a set of uniformly
slowly increasing infinitely differentiable $\Mmm$-valued functions on~$\Rn$.}\eqno{(2.8)^*}
$$
\end{props}

Our proofs of Propositions \ref{props2.2} and \ref{props2.3} are based on the Shilov inequality. In \cite[Sec.~I.5]{P}, in the proof of the prototype of Proposition \ref{props2.2}, instead of the Shilov inequality, I.~G.~Petrovski\u\i\  used \cite[Sec.~I.5, Lemma 5]{P}. We shall prove Propositions \ref{props2.2} and \ref{props2.3} according to the schemes \eqref{eq2.3}$\Rightarrow$\eqref{eq2.4}\break$\Rightarrow$\eqref{eq2.5}$\Rightarrow$\eqref{eq2.4}$\Rightarrow$\eqref{eq2.3} and \eqref{eq2.6}$\Rightarrow$\eqref{eq2.7}$\Rightarrow$\eqref{eq2.8}$\Rightarrow$\eqref{eq2.8}*$\Rightarrow$\eqref{eq2.7}$^*$$\Rightarrow$\eqref{eq2.6} where the implications \eqref{eq2.5}$\Rightarrow$\eqref{eq2.4} and \eqref{eq2.8}$\Rightarrow$\eqref{eq2.8}$^*$$\Rightarrow$\eqref{eq2.7}$^*$ are trivial.

\begin{proof}[Proof of \eqref{eq2.3}$\Leftrightarrow$\eqref{eq2.4}]
If $A(\cdot)\in\calO_M(\Rn;\Mmm)$ and \eqref{eq2.3} holds, then, by \eqref{eq2.1} and \eqref{eq2.2}, for any fixed $T\in\]0,\infty\[$ there are $C,D\in\]0,\infty\[$ and
$l\in\symN_0$ such that for every $(t,\xi)\in[0,T]\times\Rn$ one has
\begin{align*} 
\|\!\exp(tA(\xi))\|_{\Mmm}
&\le e^{t\max\hRe\sigma(A(\xi))}\bigg(1+\sum_{k=1}^{m-1}\frac{(2t)^k}{k!}\|A(\xi)\|^{k}_{\Mmm}\bigg)\\
&\le e^{TC(1+\log(1+|\xi|))}(1+2T\|A(\xi)\|_{\Mmm})^{m-1}\\
&\le D(1+|\xi|)^{TC+l(m-1)},
\end{align*}
so that \eqref{eq2.4} is satisfied. Conversely, if \eqref{eq2.4} holds, then there are
$C\in\]0,\infty\[$ and $k\in\symN_0$ such that $\|\!\exp A(\xi)\|_{\Mmm}\le C(1+|\xi|)^k$ 
for every $\xi\in\Rn$, whence, by \eqref{eq2.2},
\begin{align*}
\max\hRe\sigma(A(\xi))&=\log\rho(\exp A(\xi))\\
&\le \log\|\!\exp A(\xi)\|_{\Mmm}
\le\log C+k\log(1+|\xi|),
\end{align*}
so that \eqref{eq2.3} holds.
\end{proof}

\begin{proof}[Proof of \eqref{eq2.6}$\Rightarrow$\eqref{eq2.7}] If \eqref{eq2.6} holds, then, by \eqref{eq2.1} and \eqref{eq2.2}, for every $t\in[0,\infty\[$ and $\xi\in\Rn$ one has
\begin{align*}
\|\!\exp(tA(\xi))\|_{\Mmm}
&\le e^{s_0t}\bigg(1+\sum_{k=1}^{m-1}\frac{(2t)^k}{k!}
\|A(\xi)\|_{\Mmm}^k\bigg)\\
&\le e^{s_0t}(1+2t)^{m-1}(1+\|A(\xi)\|_{\Mmm})^{m-1}. 
\end{align*}
Furthermore, since $A(\cdot)\in\calO_M(\Rn;\Mmm)$, there are $C\in\]0,\infty\[$ and
$l\in\symN_{0}$ such that $\|A(\xi)\|_{\Mmm}\le C(1+|\xi|)^l$ for every
$\xi\in\Rn$. The above inequalities imply \eqref{eq2.7}.
\end{proof}

\begin{proof}[Proof of $\eqref{eq2.7}^*$$\Rightarrow$\eqref{eq2.6}] By \eqref{eq2.2},
$$
\max\hRe\sigma(A(\xi))=\frac1t\log\rho(\exp(tA(\xi)))\le \frac1t\log
\|\!\exp(tA(\xi))\|_{\Mmm}
$$ 
for every $t\in\]0,\infty\[$ and $\xi\in\Rn$. So, if \eqref{eq2.7}$^*$
holds, then for every $\varepsilon>0$ there are $C\in\]0,\infty\[$ and
$k\in\symN$ such that\vg
$$
\max\hRe\sigma(A(\xi))\le s_0+\varepsilon+\frac1t\log(C(1+|\xi|)^k)
$$
for every $t\in\]0,\infty\[$ and $\xi\in\Rn$, whence \eqref{eq2.6} follows.
\end{proof}

\begin{proof}[Proof of \eqref{eq2.4}$\Rightarrow$\eqref{eq2.5} and 
\eqref{eq2.7}$\Rightarrow$\eqref{eq2.8}]
The proofs of these implica\-tions\break  
are similar, and both base on the argument of I.~G.~Petrovski\u\i\  from the proof of \cite[Sec.~I.2, Lemma 2]{P}. We shall limit ourselves to \eqref{eq2.7}$\Rightarrow$\eqref{eq2.8}.

For every $\alpha\in\symN_0^n$ let
$$
U_{\alpha,t}(\xi)=(\partial/\partial\xi)^\alpha\exp(tA(\xi)).
$$
Consider the condition
$$
\textbot{there is $k_\alpha\in\symN_0$ such that for every $\varepsilon>0$ there is $C_{\alpha,\varepsilon}$ in $\]0,\infty\[$
such that whenever $(t,\xi)\in[0,\infty\[\times\Rn$, then
$$
\|U_{\alpha,t}(\xi)\|_{\Mmm}\le C_{\alpha,\varepsilon}e^{(s_0+\varepsilon)t}
(1+|\xi|)^{k_\alpha}.
$$}\eqno{(2.9)_\alpha}
$$
%\enumer
\setcounter{equation}{10}
Then \eqref{eq2.7} means that (2.9)$_0$ holds, and \eqref{eq2.8} means that (2.9)$_\alpha$ holds for every $\alpha\in\symN_0^n$. So, still assuming that
$A(\cdot)\in\calO_M(\Rn;\Mmm)$, we have to prove that (2.9)$_0$ implies (2.9)$_\alpha$ for every $\alpha\in\symN_0^n$. We proceed by induction on the length of~$\alpha$. By \eqref{eq2.7}, (2.9)$_0$ is satisfied. Suppose  that (2.9)$_\beta$ is satisfied whenever $|\beta|\le l$, and take $\alpha\in\symN_0^n$ such that
$|\alpha|=l+1$. To prove (2.9)$_\alpha$, put
$$
V_{\alpha,t}(\xi)=\sum_{\beta\le\alpha,\,|\beta|\le l}\binom{\alpha}{\beta}%{\alpha\choose \beta}
\bigg(\bigg(\frac{\partial}{\partial\xi}\bigg)^{\alpha-\beta}
A(\xi)\bigg)U_{\beta,t}(\xi).
$$
Since $A(\cdot)\in\calO_M(\Rn;\Mmm)$ and (2.9)$_\beta$ holds whenever $|\beta|\le l$, it follows that
$$
 \textbot{there is $h_\alpha\in\symN_0$  such that for every
$\varepsilon>0$ there is $D_{\alpha,\varepsilon}\in\]0,\infty\[$ 
such that whenever $(t,\xi)\in[0,\infty\[\times\Rn$, then
$$
\|V_{\alpha,t}(\xi)\|_{\Mmm}\le D_{\alpha,\varepsilon} e^{(s_0+\varepsilon)t}
(1+|\xi|)^{h_\alpha}.
$$}\eqno{(2.10)_\alpha}
$$
One has
\begin{align*}
\frac{\partial}{\partial t}U_{\alpha,t}(\xi)
&=\frac{\partial}{\partial t}\bigg(\frac{\partial}{\partial\xi}\bigg)^{\alpha}
\exp(tA(\xi))
=\bigg(\frac{\partial}{\partial\xi}\bigg)^{\alpha}[A(\xi)\exp(tA(\xi))]\\
&=A(\xi)U_{\alpha,t}(\xi)+V_{\alpha,t}(\xi)
\end{align*}
 and $U_{\alpha,0}(\xi)=0$ because $|\alpha|=l+1\ge 1$. Hence
\begin{equation}
\label{eq2.11}
U_{\alpha,t}(\xi)=\int_0^t[\exp((t-\tau)A(\xi))]V_{\alpha,t}(\xi)\,d\tau.
\end{equation}
From (2.9)$_0$, (2.10)$_{\alpha}$ and \eqref{eq2.11} it follows that 
\begin{multline*}
\|U_{\alpha,t}(\xi)\|_{\Mmm}\\
\begin{aligned}
&\le \int_0^t C_{0,\varepsilon/2} 
e^{(s_0+\varepsilon/2)(t-\tau)}(1+|\xi|)^{k_0}D_{\alpha,\varepsilon/2} e^{(s_0+\varepsilon/2)\tau}(1+|\xi|)^{h_\alpha}\,d\tau\\
&=C_{0,\varepsilon/2} D_{\alpha,\varepsilon/2}te^{(s_0+\varepsilon/2)t}(1+|\xi|)^{k_0+h_\alpha}
\le \tilde C_{\alpha,\varepsilon}e^{(s_0+\varepsilon)t}(1+|\xi|)^{k_\alpha}
\end{aligned}
\end{multline*}
for $k_\alpha=k_0+h_\alpha$ and $\tilde C_{\alpha,\varepsilon}=C_{0,\varepsilon/2}
D_{\alpha,\varepsilon/2}\max_{t\in[0,\infty\[}t e^{-(\varepsilon/2)t}$.
\end{proof}

\section{Proof of Theorem 1}\label{sec3}
\looseness2
{\it Necessity of} \eqref{eq1.11}. Suppose that $(S_{t})_{t\ge0}\subset\calO_{C}^{\prime}(\Rn;\Mmm)$ is an i.d.c.s. with generating distribution $G\in\calO_{C}^{\prime}(\Rn;\Mmm)$. Let $A=\scrF G$. Then $A,\scrF S_{t}\in
\calO_{M}(\Rn;\Mmm)$ and $(\scrF S_{t})(\xi)=\exp(tA(\xi))$
for every $t\in[0,\infty\[$ and $\xi\in\Rn$. Since the mapping $[0,\infty\[\ni t\mapsto
[\exp(tA(\cdot))]\,\cdot=(\scrF S_{t})\,\cdot\in L_{b}(\calS(\Rn;\symC^{m});\calS(\Rn;\symC^{m}))$
is continuous, the Banach--Steinhaus theorem implies that whenever $T\in\]0,\infty\[$, then the set of multiplication operators $\{[\exp(tA(\cdot))]\,\cdot:t\in[0,T]\}$ is an equicontinuous subset of $L(\calS(\Rn;\symC^{m});\allowbreak \calS(\Rn;\symC^{m}))$. By \cite[Theorem 3.1]{K3}, this is equivalent to \eqref{eq2.5}. By Proposition 2.2, \eqref{eq2.5} is equivalent to \eqref{eq2.3}. Since
$A=\scrF G$, \eqref{eq2.3} is nothing but~\eqref{eq1.11}.
\ods

{\it Sufficiency of} \eqref{eq1.11}. Suppose that
$G\in\calO_{C}^{\prime}(\Rn;\Mmm)$  satisfies \eqref{eq1.11}. Let $A=\scrF G$. Then
$A\in\calO_{M}(\Rn;\Mmm)$, and $A$ satisfies \eqref{eq2.3}. Hence, by Proposition \ref{props2.2} and 
\cite[Theorem 3.1]{K3}, whenever $T\in\]0,\infty\[$, then
$\{[\exp(tA(\cdot))]\,\cdot:t\in[0,T]\}$ is an equicontinuous subset of
$L(\calS(\Rn;\symC^{m});\allowbreak\calS(\Rn;\symC^{m}))$. By the theorem on differentiating a solution of an ODE with respect to a parameter \cite[Sec.~V.4, Corollary 4.1]{Ha}, the mapping
$\symR^{1+n}\ni(t,\xi)\mapsto\exp(tA(\xi))\in\Mmm$ is infinitely differentiable, and hence, by 
\cite[Theorem 3.2]{K3}, so is $[0,\infty\[\ni t\mapsto[\exp(tA(\cdot))]\,\cdot\in
L_{b}(\calS(\Rn);\calS(\Rn))$, and its right derivative at zero (computed in the topology of 
$L_{b}(\calS(\Rn;\symC^{m});\allowbreak\calS(\Rn;\symC^{m}))$) is
$A\,\cdot\in L(\calS(\Rn);\calS(\Rn))$. It follows that
$G\,*=(\scrF^{-1}A)\,*=\scrF^{-1}\circ(A\,\cdot)\circ\scrF\in
L(\calS(\Rn;\symC^{m});\calS(\Rn;\symC^{m}))$ is the infinitesimal generator of the infinitely differentiable operator semigroup 
$([\scrF^{-1}\exp(tA(\cdot))]\,*)_{t\ge0}=
(\scrF^{-1}\circ[\exp(tA(\cdot))]\,\cdot)\circ\scrF)_{t\ge0}\subset
L_{b}(\calS(\Rn;\symC^{m});\calS(\Rn;\symC^{m}))$. Consequently, $G=\scrF^{-1}A$ is the generating distribution of the i.d.c.s. $(\scrF^{-1}\exp(tA(\cdot)))_{t\ge0}\subset\calO_{C}^{\prime}(\Rn;\Mmm)$.

\section{Proof of Theorem 2}\label{sec4}

Let $(S_t)_{t\ge0}\subset \calO_C^\prime(\Rn;\Mmm)$ be an i.d.c.s. with generating distribution
$G\in\calO_C^\prime(\Rn;\Mmm)$. 
Put $A=\scrF G$. Then $A,\scrF S_t\in\calO_M (\Rn;\Mmm)$, 
condition (b) from the Lemma from Sec.~1.4 is satisfied, and
$(\exp(tA(\cdot)))\,\cdot=\scrF\circ(S_t\,*)\allowbreak\circ\scrF^{-1}$ for every $t\in[0,\infty\[$.
Since $\scrF,\scrF^{-1}\!\in\! L(\calS(\Rn);
\calS(\Rn))$, 
for $\omega((S_t)_{t\ge0})$ defined by \eqref{eq1.13}
one has 
\begin{multline*}
\indent
\omega((S_t)_{t\ge 0})=\inf\{\omega\in\symR:
\{[e^{-\omega t}\exp(tA(\cdot))]\cdot:t\in[0,\infty\[\}
\mbox{ is}\\
\mbox{an equicontinuous subset of }
 L_b(\calS(\Rn;\symC^m);\calS(\Rn;\symC^m))\}.
\end{multline*}	 
From \cite[Theorem 3.1]{K3} it follows that whenever $s_0\in\symR$, then 
\begin{equation}\label{eq4.1}
\omega((S_t)_{t\ge0})<s_0+\varepsilon\ \quad\mbox{for every }\varepsilon>0
\end{equation}
if and only if \eqref{eq2.8}$^*$ holds. Hence, by Proposition \ref{props2.3}, the condition \eqref{eq4.1} is equivalent to \eqref{eq2.6}. This implies that the growth bound of the i.d.c.s.
$(S_t)_{t\ge0}\subset\calO_C^\prime(\Rn;\Mmm)$ is equal to the spectral bound of~$G$, where both these quantities may well be infinite. 

\section{Condition (\ref{eq1.22}) implies G{\aa}rding hyperbolicity}\label{sec5}

Let $(S_t)_{t\ge0}\subset\calO_C^\prime(\Rn;\Mmm)$ be an i.d.c.s. with generating distribution $G=\calG(\partial_1,\ldots,\partial_n)\otimes\delta$, so that the condition \eqref{eq1.14} is satisfied. Then
$\scrF S_t\in\calO_M(\Rn;\Mmm)$ and $(\scrF S_t)(\xi)=\exp(t\calG(i\xi))$ for every $t\in[0,\infty\[$ and $\xi\in\Rn$. Suppose that \eqref{eq1.22} holds, i.e.
$S_{t_0}\in\calE^\prime(\Rn;\Mmm)$
for some $t_0\in\]0,\infty\[$. Then, by the Paley--Wiener--Schwartz theorem, i.e. by \cite[Theorem 7.3.1]{H} or \cite[Theorem 8.57]{K-R}, there are
$C,k,l\in\]0,\infty\[$ such that whenever $\zeta\in\symC^n$, then 
\begin{equation}
\label{eq5.1}
\|\!\exp(t_0\calG(i\zeta))\|_{\Mmm}
=\|(\scrF S_{t_0})(\zeta)\|_{\Mmm}
\le C(1+|\zeta|)^l e^{k\hIm \zeta}.
\end{equation}
For every $\zeta\in\symC^n$ put
\begin{equation*}
\Lambda(\zeta)=\max\hRe\sigma(\calG(i\zeta)).
\end{equation*}
 Then
\begin{equation*}
\Lambda(\zeta)=\max\{\hRe\lambda:\lambda\in\symC,\, 
P (\lambda,\zeta_1,\ldots,\zeta_n)=0\}
\end{equation*}
where
\begin{align*}
P (\lambda,\zeta_1,\ldots,\zeta_n)
={}&\det(\lambda\symjed_{\mm}-\calG(\zeta_1,\ldots,\zeta_n))\\
={}&\lambda^m+Q_{m-1}(\zeta_1,\ldots,\zeta_n)\lambda^{m-1}\\
&{}+\cdots+
Q_{1}(\zeta_1,\ldots,\zeta_n)\lambda+Q_0(\zeta_1,\ldots,\zeta_n).
\end{align*}
Let
$$
p_0=\inf\{p\in\]0,\infty\[:\sup_{\zeta\in\symC^n}(1+|\zeta|)^{-p}
\Lambda(\zeta)<\infty\}.
$$
By \eqref{eq2.2} and \eqref{eq5.1} there is $K\in \]0,\infty\[$ such that
$$
\Lambda(\zeta)\le t_0^{-1}\log\|\!\exp(t_0\calG(i\zeta))\|_{\Mmm}
\le K(1+|\zeta|)
$$
for every $\zeta\in\symC^n$. Consequently,
\begin{equation}\label{eq5.2}
p_0\le 1.
\end{equation} 
By the Gelfand--Shilov theorem on the reduced order \cite[Sec.~II.6.2]{G-S}, \cite[Sec.~7.2, Theorem~4]{F},
$$
p_0=\max_{k=0,\ldots,m-1}(m-k)^{-1}\deg Q_k,
$$
so that, by \eqref{eq5.2}, $\deg Q_k\le m-k$ for every $k=0,\ldots,m-1$, and hence $\deg P=m$, proving \eqref{eq1.23}.

\section{G{\aa}rding hyperbolicity implies Ehrenpreis hyperbolicity}\label{sec6}

Suppose that $(S_t)_{t\ge0}\subset\calO_C^\prime(\Rn;\Mmm)$ is an i.d.c.s. with generating distribution $G=\calG(\partial_1,\ldots,\partial_n)\otimes\delta$. Let $P(\lambda,\zeta_1,\ldots,\zeta_n)
=\det(\lambda\symjed_{\mm}-\calG(\zeta_1,\allowbreak\ldots,\zeta_n))$. Then, by Theorem 3, \eqref{eq1.14}$^\prime$ holds, i.e. $\sup\{\hRe\lambda:\lambda\in\symC$ and there is
$(\xi_1,\ldots,\xi_n)\allowbreak\in\Rn$ such that
$P(\lambda,i\xi_1,\ldots,i\xi_n\}=0\}=s_0(\calG)<\infty$. Suppose moreover that \eqref{eq1.23} holds, i.e. $\deg P(\lambda,\zeta_{1},\ldots,\zeta_{n})=m$.

By \cite[Theorem 12.4.2 and Lemma 8.7.3]{H}, the above properties of\break
$P(\lambda,\zeta_1,\ldots,\zeta_n)$ imply that $\Gamma$ defined in our Sec.~1.5 is an open convex cone with vertex at zero. 
From the definition of $\Gamma$ it follows that $\Gamma$ contains the open halfline
$\{(t,0,\ldots,0)\in\symR^{1+n}:t>0\}$. 
From \cite[Theorem 12.4.4]{H}~$^{*)}$\footnote{$^{*)}$ One could also use
\cite[Lemma~2.6]{G}, but in \cite{G} the open convex cone
$\Gamma$ has a definition equivalent to but formally different
from ours, which is taken from~\cite{H}.} it follows that
\begin{equation}
\textbot{whenever $(\nu_0,\nu_1,\ldots,\nu_n)\in\Gamma$,
$(\xi,\ldots,\xi_n)\in\Rn$,
${\lambda,\mu\in\symC}$,
$\hRe\lambda>s_0(\calG)$ and $\hRe\mu\ge0$, then 
$$P(\lambda+\mu\nu_0,i\xi_1+\mu\nu_1,\ldots,i\xi_n+\mu\nu_n)\ne0.$$}\label{eq6.1}
\end{equation}
Fix $r>0$ so large that
$$
K_r:=\{(\nu_0,\nu_1,\ldots,\nu_n)\in\symR^{1+n}:\nu_0\ge r,\,\nu_1^2+\cdots+
\nu_n^2\le1\}\subset\Gamma.
$$
Let $(\xi_1,\ldots,\xi_n),(\eta_1,\ldots,\eta_n)\in\Rn,\,\mu=1+|\eta|
=1+(\eta_1^2+\cdots+\eta_n^2)^{1/2},(\nu_0,\nu_1,\allowbreak\ldots,\nu_n)
=(r,\eta_1/(1+|\eta|),\ldots,\eta_{n}/(1+|\eta|))$.
Then $(\nu_0,\ldots,\nu_n)\in K_r\subset\Gamma$, and if $\lambda\in\symC$ and
$$
\hRe\lambda> s_0(\calG)+(1+|\eta|)r,
$$
then, by \eqref{eq6.1},
\begin{multline*}
P(\lambda,i\xi_1+\eta_1,\ldots,i\xi_n+\eta_n)\\
=P((\lambda-(1+|\eta|)r)+\mu\nu_0,i\xi_1+\mu\nu_1,\ldots,i\xi_n+\mu\nu_n)\ne0
\end{multline*}
because $\hRe(\lambda-(1+|\eta|)r)>s_0(\calG)$. It follows that 
$$
 \textbot{whenever $(\lambda,\zeta_1,\ldots,\zeta_n)\in\symC^{1+n}$
and $P(\lambda,\zeta_1,\ldots,\zeta_n)=0$, then 
$$
\hRe\lambda\allowbreak\le s_0(\calG)+r+r((\hRe\zeta_1)^2+\cdots+(\hRe\zeta_n)^2)^{1/2}.
$$}\eqno{(6.2)_+}
$$
By \cite[Lemma 2.2]{G} or \cite[Theorem 12.4.1]{H}, if the polynomial
$P(\lambda,\zeta_1,\ldots,\zeta_n)$ satisfies \eqref{eq1.14}$^\prime$ and \eqref{eq1.23}, then so does $P(-\lambda,\zeta_1,\ldots,\zeta_n)$. Since (6.2)$_+$ is a consequence of the properties \eqref{eq1.14}$^\prime$ and \eqref{eq1.23} of $P(\lambda,\zeta_1,\ldots,\zeta_n)$, it follows that the properties \eqref{eq1.14}$^\prime$ and \eqref{eq1.23} of $P(-\lambda,\zeta_1,\ldots,\zeta_n)$ imply that there is $r^\prime>0$ such that
$$
\textbot{whenever $(\lambda,\zeta_1,\ldots,\zeta_n)\in\symC^{1+n}$ and 
$P(\lambda,\zeta_1,\ldots,\zeta_n)=0$, then
$$
-\hRe\lambda\allowbreak\le s_0(-\calG)+r^\prime+r^\prime((\hRe\zeta_1)^2+\cdots+(\hRe\zeta_n)^2)^{1/2}.
$$}\eqno{(6.2)_-}
$$
Together $(6.2)_+$ and $(6.2)_-$ mean that \eqref{eq1.27} is satisfied, i.e. the matricial PDO \eqref{eq1.26} is hyperbolic in the sense of Ehrenpreis with respect to the coordinate~$t$.

\section{\!\!\!The Ehrenpreis hyperbolicity implies (\ref{eq1.24})}\label{sec7}

Suppose that the system \eqref{eq1.26} is hyperbolic in the sense of Ehrenpreis with respect to the coordinate~$t$. This means that whenever $(\lambda,\zeta_1,\ldots,\zeta_n)\in\symC^{1+n}$ and 
$$
P(\lambda,\zeta_1,\ldots,\zeta_n)
=\det(\lambda\symjed_{\mm}-\calG(\zeta_1,\ldots,\zeta_n))=0,
$$
then
$$
|\!\hRe\lambda|\le C(1+((\hRe\zeta_1)^2+\cdots+(\hRe\zeta_n)^2)^{1/2})
$$
for some $C\in\]0,\infty\[$ independent of $(\lambda,\zeta_1,\ldots,\zeta_n)$. Since
$$
\sigma(\calG(i\zeta))=\{\lambda\in\symC:P(\lambda,i\zeta,\ldots, i\zeta)=0\},
$$
it follows that whenever $(\zeta_1,\ldots,\zeta_n)\in\symC^n$, then
\begin{equation}\label{eq7.1}
\max|\!\hRe\sigma(\calG(i\zeta))|\le C(1+((\hIm\zeta_1)^2+\cdots+
(\hIm\zeta_n)^2)^{1/2}).
\end{equation}
By \eqref{eq2.1} and \eqref{eq2.2}, this implies that
\begin{align}
\|\!\exp(t\calG(i\zeta))\|_{\Mmm}
&\le e^{C|t|}\bigg(1+\sum_{k=1}^{m-1}\frac{(2|t|)^k}{k!}
\|\calG(i\zeta)\|_{\Mmm}^k\bigg) e^{C|t|\,|\!\hIm\zeta|}\notag\\
&\le e^{C|t|}(1+2|t|)^{m-1} D(1+|\zeta|)^{(m-1)d}e^{C|t|\,|\!\hIm\zeta|}\label{eq7.2},
\end{align}	 
for every $(t,\zeta)\in\symR\times\symC$ where $C,D\in\]0,\infty\[$ are independent of $(t,\zeta)$, and  $d\in\symN_0$ is the maximum of the orders of the scalar PDO which are the entries of $\calG(\partial_1,\ldots,\partial_n)$. By the Paley--Wiener--Schwartz theorem, i.e. by \cite[Theorem 7.3.1]{H}, \eqref{eq7.2} implies that there is a one-parameter convolution group $(\tilde S_t)_{t\in\symR}\subset\calE^\prime(\Rn;\Mmm)$ such that
\begin{equation}\label{eq7.3}
(\scrF\tilde S_t)(\zeta)=\exp(t\calG(i\zeta))\ \quad\mbox{for every }
(t,\zeta)\in\symR\times\symC^n
\end{equation}
and
$$
\max\{|x|:x\in\supp \tilde S_t\}\le C|t|\ \quad\mbox{for every }
t\in\symR.
$$
The convolution group $(\tilde S_t)_{t\in\symR}$ is an extension of the i.d.c.s. $(S_t)_{t\ge0}\subset\calO_C^\prime(\Rn;\Mmm)$ with generating distribution
$\calG(\partial_1,\ldots,\partial_n)\otimes\delta$ which exists by Theorem \ref{thm3} because \eqref{eq7.1}$\Rightarrow$\eqref{eq1.14}. Furthermore, by \eqref{eq7.2}, one has
\begin{multline}
\label{eq7.4}
\|\calG(i\zeta)^k\exp(t\calG (i\zeta))\|_{\Mmm}\\
\le e^{C|t|}(1+2|t|)^{m-1}D_k(1+|\zeta|)^{(m+k-1)d}e^{C|t|\,|\!\hIm\zeta|}
\end{multline}	
for every $(t,\zeta)\in\symR\times\symC^n$ and $k\in\symN_0$. By the theorem on differentiating a solution of an ODE with respect to a parameter (\cite[Sec.~V.4, Corollary 4.1]{Ha}), the mapping $\symR\times\symC^n\ni(t,\zeta)\mapsto\exp(t\calG(i\zeta))\in\Mmm$
 is infinitely differentiable. Since $(\partial/\partial t)^k\exp(t\calG(i\zeta))=\calG(i\zeta)^k\exp(t\calG(i\zeta))$, from \eqref{eq7.4} and \cite[Sec.~V.5, Lemma 5.17]{E2} it follows that the mapping $\symR\ni t\mapsto\tilde S_t\in\calE^\prime(\Rn;\Mmm)$ is infinitely differentiable in the topology of $\calE^\prime(\Rn;\Mmm)$.

\section{Application to the Cauchy problem}

\subsection{Well posedness spaces}
Let $\calG(\partial_1,\ldots,\partial_n)$ be an $\mm$ matrix whose entries are PDOs on  $\Rn$ with constant complex coefficients. Suppose that 
\begin{equation}
\sup\{\hRe\lambda:\lambda\in\sigma(\calG(i\xi)),\,\xi\in\Rn\}<\infty. \tag{ii}
\end{equation}
{\spaceskip.33em plus.22em minus.17em Then, by Theorem 3, there is a unique infinitely differentiable convolution semigroup $(S_t)_{t\ge0}\subset\calO_C^\prime(\Rn;\Mmm)$
with generating distribution $\calG(\partial_1,\ldots,\partial_n)\otimes\delta$. Suppose moreover that}
%\iiinumer
$$
 \textbot{$E$ is a sequentially complete l.c.v.s. continuously imbedded in\break
${\calS^\prime(\Rn;\symC^m)}$ such
that $(S_t\,*)E\subset E$ for every $t\in[0,\infty\[$,\ and the mapping
$[0,\infty\[\times E\ni(t,u)\mapsto S_t*u\in E$ is separately continuous.}\eqno{{\rm(iii)}_{S_t,E}}
$$
Define the operator $\calG_E$ from $E$ into $E$ by the conditions
\begin{align*}
D(\calG_E)&=\{u\in E:\calG(\partial_1,\ldots,\partial_n)u\in E\}, \\
\calG_Eu&=\calG(\partial_1,\ldots,\partial_n)u\ \quad\mbox{for }u\in D(\calG_E).
\end{align*}

\begin{thm}\label{thm5}
{\it Suppose that conditions {\rm(ii)} and {\rm(iii)}$_{S_t,E}$ are satisfied. Then for every $k=1,2,\ldots,\infty$ the Cauchy problem
\begin{equation}
\frac{d}{dt}u(t)=\calG(\partial_1,\ldots,\partial_n)u(t)\ \quad
\mbox{for }t\in[0,\infty\[,\ \quad u(0)=u_0, \tag{iv}
\end{equation}
with given $u_0\in D(\calG_E^k)$ has a solution $u(\cdot)\in C^k([0,\infty\[;E)$ which is unique in the class $C^1([0,\infty\[;\calS^\prime(\Rn;\symC^m))$. This solution is given by the formula}
\begin{equation}
u(t)=S_t*u_0\ \quad\mbox{\it for }t\in[0,\infty\[.\tag{v}
\end{equation}
\end{thm}

Thanks to Theorem \ref{thm5} it is legitimate to call $E$ the {\it well posedness space} for the Cauchy problem (iv) if conditions (ii) and ${\rm(iii)}_{S_t,E}$ are satisfied. Theorem \ref{thm5} confirms the observation of L.~H\"ormander \cite[notes at the end of Chapter 12]{H} that the Petrovski\u\i\ condition (ii) is related to well posedness of the Cauchy problem for PDOs with constant coefficients in L.~Schwartz spaces $\calS$ and~$\calS^\prime$.
\ods

\begin{xrem} Let $Z(\symC^n;\symC^m)$ be the space of $\symC^m$-valued functions holomorphic on $\symC^n$ such that $\varphi\in Z(\symC^n;\symC^m)$ if and only if there is $a=a(\varphi)\in\]0,\infty\[$ such that
$\sup_{z\in\symC^n}(1+\|z\|)^k e^{-a\|\!\hIm z\|}\|\varphi(z)\|<\infty$
for every $k\in\symN$. Let $Z(\symC^n;\symC^m)|_{\Rn}$ be the set of restrictions to $\Rn$ of functions in $ Z(\symC^n;\symC^m)$. By the Paley--Wiener theorem
(\cite[Theorem 7.3.1]{H}, \cite[Theorem 8.51]{K-R}), $Z(\symC^n;\symC^m)|_{\Rn}=\scrF\calD(\Rn;\symC^{m})$, so that $Z(\symC^n;\symC^m)|_{\Rn}$ is a dense subset of $\calS(\Rn;\symC^m)=\scrF\calS(\Rn;\symC^m)$. {\it If $E=Z(\symC^n;\symC^m)|_{\Rn}$, then the Cauchy problem {\rm(iv)} is well posed for every $\calG(\partial_1,\ldots,\partial_n)$ independently of whether {\rm(ii)} holds or not.} 
Indeed, if $E=Z(\symC^n;\symC^m)|_{\Rn}$, then (instead of appealing to Sec.~\ref{sec2} which enables the use of the Lemma from Sec.~1.4) in order to conclude that the Cauchy problem (iv) is well posed it is sufficient to observe that the mapping
$
\symR\ni t\mapsto\[\exp t\calG(i\,\cdot,\ldots,i\,\cdot)]\,\cdot\in L_b(\calD(\Rn;\symC^m);\calD(\Rn;\symC^m))
$
is infinitely differentiable because so is the mapping
$\symR^{1+n}\ni(t,\xi_1,\ldots,\xi_n)\mapsto\exp (t\calG(i\xi_1,\ldots,i\xi_n))
\in\Mmm$.
\end{xrem}

\subsection{Examples of well posedness spaces}

Examples of spaces $E$ satisfying ${\rm(iii)}_{S_t,E}$ for each
$\calG(\partial_1,\ldots,\partial_n)$ satisfying (ii) include:
\begin{aitem}
\item[(a)] the spaces of infinitely differentiable functions $\calS(\Rn;\Cm)$ and\break
$\calD_{L^p}(\Rn;\Cm)\!=\!\{u\!\in\! C^\infty(\Rn;\Cm):\partial^\alpha u\!\in\! L^p(\Rn;\Cm)$ 
for every $\alpha\in\nobreak\symN_0^n\}$, $p\in[1,\infty]$,
\item[(b)] the spaces of distributions $\calS'(\Rn;\Cm)$, $\calO'_C(\Rn;\Cm)$ and
$$
\calD'_{L^q}(\Rn;\Cm)=(\calD_{L^p})'(\Rn;\Cm),\quad
q\in\]1,\infty],\quad
p=q/(q-1).
$$
\end{aitem}
Examples of spaces $E$ depending on $\calG(\partial_1,\ldots,\partial_n)$ such that the Cauchy problem (iv) is well posed whenever $\calG(\partial_1,\ldots,\partial_n)$ satisfies (ii) include:
\begin{aitem}
\item[(c)] the T.~Ushijima space\vadjust{\vskip-4pt}
\begin{multline*}
U_{\calG}(\Rn;\Cm)=\{u\in L^2(\Rn;\Cm):(\calG(\partial_1,\ldots,\partial_n))^k
u\in L^2(\Rn;\Cm)\\				
\mbox{for every }k=1,2,\ldots\}
\end{multline*}\vskip-\lastskip
occurring in \cite[Theorem 10.1]{U},
\item[(d)] the Banach spaces $\calB_{\calN,p}$ of G.~Birkhoff \cite{B},
\item[(e)] the Hilbert spaces $\calL_{B}$ of S.~D.~Eidelman and S.~G.~Krein discussed in 
\cite[Sec.~I.8.2]{K}.
\end{aitem}
In the cases (c)--(e) the well posedness of the Cauchy problem (iv) follows directly from the results of \cite{U}, \cite{B} and \cite{K} without reference to
$\calO'_C(\Rn;\Cm)$ and (v). Notice that in \cite{P}, \cite{U} and \cite{K2} it is proved that if $E$ is equal to either of the spaces
$\calD_{L^\infty}(\Rn;\Cm)$, $\calD_{L^2}(\Rn;\Cm)$ or
$U_{\calG}(\Rn;\Cm)$, then (ii) is necessary for well posedness of the Cauchy problem (iv). (The arguments from \cite{P} and \cite{U} are quoted in \cite{K2}.)

From among the spaces $\calD_{L^p}(\Rn;\Cm)$, $p\in[1,\infty]$, the largest
one is $\calD_{L^\infty}(\Rn;\Cm)$ whose dual is not a space of distributions and does not occur in~(b). The fact that if (ii) holds, then the Cauchy problem (iv) 
is well posed for $E=\calD_{L^\infty}(\Rn;\Cm)$, was first proved by I.~G.~Petrovski\u\i\
\cite{P} in 1938. From among the spaces $\calD'_{L^q}(\Rn;\Cm)$, $q\in\]1,\infty]$, 
the largest one is $\calD'_{L^\infty}(\Rn;\Cm)$, i.e. the space of $\Cm$-valued 
distributions on $\Rn$ bounded in the sense of L.~Schwartz. 
An alternative notation for $\calD'_{L^\infty}$ is~$\calB'$.

The space $\calO_M(\Rn;\Cm)$ cannot be included in (a) because
$\calO_M(\Rn;\mathbb{C}^{1})$ is not a well posedness space when
$\calG(\partial_1,\ldots,\partial_n)=i\Delta$. Indeed $i\Delta$ satisfies (ii) and, 
in accordance with Example 1 of Sec.~1.4, the i.d.c.s. whose generating distribution is
$i\Delta \delta$ extends to a one-parameter group
$(S_t)_{t\in\symR}\subset \calO'_C(\Rn)$ such that
$S_t\in \calO'_C(\Rn)\cap\calO_M(\Rn)$ for every $t\in\symR\setminus\{0\}$. 
Fix $t_0\in\]0,\infty\[$. Then, by Theorem~5, the Cauchy problem
$$
\frac{d}{dt}u(t)=i\Delta u(t)\quad\mbox{for }t\in[0,\infty\[,\quad u(0)=S_{-t_0},
\eqno{{\rm(iii)}_0}$$
has in the class $C^1([0,\infty\[;\calS'(\Rn))$ a unique solution. 
Since this unique solution is given by the formula
$u(t)=S_t*S_{-t_0}$ it follows that $u(0)=S_{-t_0}\in\calO_M(\Rn)$  and 
$u(t_0)=\delta\not\in\calO_M(\Rn)$. Consequently, the Cauchy problem (iii)$_0$  has no solution in the class $C^1([0,\infty\[;\calO_M(\Rn))$.

\subsection{Well posedness of the spaces  $\calD_{L^p}(\Rn;\Cm)$}

We shall use the following

\begin{xlem} Let $(S_t)_{t\ge 0}\subset \calO_C(\Rn;\Mmm)$ be an i.d.c.s. 
with generating distribution $\calG(\partial_1,\ldots,\partial_n)\otimes\delta$ satisfying {\rm(ii)}. Then there are $j_0\in\symN$ and a continuous mapping $[0,\infty\[\ni t\mapsto f_t\in L^1(\Rn;\Mmm)$
having the three properties:
\begin{aitem}
\item[\rm(a)] $f_t\in L^1(\Rn;\Mmm)\cap\calO'_C(\Rn;\Mmm)$ for every $t\in[0,\infty\[$,
\item[\rm(b)] $S_t=(1-\Delta)^{j_0}f_t$ for every $t\in[0,\infty\[$ where the right side is understood in the sense of $\calS'(\Rn;\Mmm)$,
\item[\rm(c)] $\sup_{t\in[0,\infty\[} e^{-(s_0(\calG)+\varepsilon)t}\|f_t\|_{L^1(\Rn;\Mmm)}<\infty$ for every $\varepsilon>0$.
\end{aitem}
\end{xlem}

Before proving the lemma let us show how it implies that $\calD_{L^p}(\Rn;\Cm)$,
$p\in[1,\infty]$, are well posedness spaces. Let
$(S_t)_{t\ge0}\subset\calO'_C(\Rn;\Mmm)$
be an i.d.c.s. whose generating distribution
$\calG(\partial_1,\ldots,\partial_n)\otimes\delta$ satisfies (ii). 
Whenever $t\in[0,\infty\[$ and $u\in \calD_{L^p}(\Rn;\Cm)$ then\vadjust{\vskip-2pt}
$$
S_t*u=((1-\Delta)^{j_0}f_t)*u=f_t*((1-\Delta)^{j_0}u)\vadjust{\vskip-2pt}
$$
because $S_t,f_t\in\calO'_C(\Rn;\Mmm)$ and $u\in \calS'(\Rn;\Cm)$. The continuity of the mapping
$[0,\infty\[\ni t\mapsto f_t\in L^1(\Rn;\Mmm)$ implies the separate continuity of the mapping\vadjust{\vskip-2pt}
$$
[0,\infty\[\times \calD_{L^p}(\Rn;\Cm)\ni(t,u)\mapsto S_t*u
=f_t*((1-\Delta)^{j_0}u)\in\calD_{L^p}(\Rn;\Cm).\vadjust{\vskip-2pt}
$$
Finally, $\calD_{L^p}(\Rn;\Cm)$ is sequentially complete, and it is continuously imbedded in $\calS'(\Rn;\Cm)$. Consequently, the condition (iii)$_{S_t,E}$ is satisfied for
$E=\calD_{L^p}(\Rn;\Cm)$.

\begin{proof}[Proof of the Lemma.] By the estimation \eqref{eq2.8} from Proposition \ref{props2.3}, and by the statement (2.8) from \cite{K3}, for every fixed
$j_0\in\symN$ and $t\in[0,\infty\[$ the function
$g_t:\Rn\ni\xi\mapsto(1+|\xi|^2)^{-j_0}\exp(t\calG(i\xi))\in\Mmm$
belongs to $\calO_M(\Rn;\Mmm$). Let $f_t=\scrF^{-1}g_t$. 
Then $f_t=\calO'_C(\Rn;\Mmm)$ and $(1-\Delta)^{j_0}f_t=\scrF^{-1}(\exp(t\calG(i\,\cdot)))
=S_t$. The Lemma follows once we prove that if $j_0$ is sufficiently large, then each distribution $f_t\in\calO'_C(\Rn;\Mmm)$, $t\in[0,\infty\[$, is represented by a function
 belonging to $L^1(\Rn;\Mmm)$ such that the mapping
$[0,\infty\[\ni t\mapsto f_t\in L^1(\Rn;\Mmm)$ is locally lipschitzian and satisfies~(c).
 We shall base on the fact that
\begin{equation}\label{eq8.1d}
\textbot{if $T\in\calS'(\Rn;\Mmm)$  and $(1-\Delta)^{[n/2]+1}\hat T\in L^1(\Rn;\Mmm)$,\hfill\break 
then $T\in L^1(\Rn;\Mmm)$ and\vadjust{\vskip-2pt} 
$$
\|T\|_{L^1(\Rn;\Mmm)}\le C\|(1-\Delta)^{[n/2]+1}\hat T\|_{L^1(\Rn;\Mmm)}\vadjust{\vskip-2pt}
$$}
\end{equation}
where $C\in\]0,\infty\[$ depends only on~$n$. To prove \eqref{eq8.1d} it is sufficient to note that one has dense imbeddings $\calS(\Rn;\Mmm)\subset L^1(\Rn;\Mmm)\subset 
\calS'(\Rn;\Mmm)$  and if $\varphi\in\calS(\Rn;\Mmm)$, then
\begin{align*}
\|\varphi\|_{L^1(\Rn;\Mmm)}&\le C\sup_{x\in\Rn}(1+|x|^2)^{[n/2]+1}\|\varphi(x)\|_{\Mmm}\\
&\le C\|(1-\Delta)^{[n/2]+1}\hat\varphi\|_{L^1(\Rn;\Mmm)}.
\end{align*}

In order to prove the Lemma, we shall apply \eqref{eq8.1d} to
$T=(d/dt)^l f_t$ where $l=0,1$ and $t\in[0,\infty\[$. For this $T$ one has
$$ 
((1-\Delta)^{[n/2]+1}\hat T)(\xi)=(1-\Delta_{\xi})^{[n/2]+1}
((1+|\xi|^{2})^{-j_{0}}(\calG(i\xi))^l\exp(t \calG(i\xi))),
$$
so that $(1-\Delta)^{[n/2]+1}\hat T\in\calO_M(\Rn;\Mmm)$, again by \eqref{eq2.8} from Proposition \ref{props2.3} and (2.8) from \cite{K3}. 
In order to show that if $j_0$ is sufficiently large, then
$(1-\Delta)^{[n/2]+1}\hat T\in  L^1(\Rn;\Mmm)$,
it is sufficient to prove that whenever $j_0\in\symN$ is sufficiently large and
$\kappa\in\symN^n_0$ is a multiindex
of length $|\kappa|\le n+2$, then the $\Mmm$-valued function 
$$
\xi\mapsto\bigg(\frac{\partial}{\partial \xi}\bigg)^{\kappa}
[(1+|\xi|^2)^{-j_0}(\calG(i\xi))^l\exp(t\calG(i\xi))]
$$
is integrable on~$\Rn$.

The Leibniz formula, the estimation \eqref{eq2.8} from Proposition \ref{props2.3}, and the statement (2.8) from \cite{K3} imply that for every $\varepsilon>0$
there is $D_\varepsilon\in\]0,\infty\[$ such that whenever $|\kappa|\le n+2$, then
\begin{multline*} 
 \sum_{|\kappa|\le n+2} \int_{\Rn}\bigg\|\bigg(\frac{\partial}{\partial\xi}\bigg)^{\kappa}
[(1+|\xi|^2)^{-j_0}(\calG(i\xi))^l\exp(t\calG(i\xi))]\bigg\|_{\Mmm}\,d\xi\\
\hskip24pt\le \sum_{|\alpha|+|\beta|+|\gamma|\le n+2}\frac{(\alpha+\beta+\gamma)!}{\alpha!\beta!\gamma!}
\int_{\Rn}\bigg|\bigg(\frac{\partial}{\partial\xi}\bigg)^{\gamma}(1+|\xi|^2)^{-j_0}\bigg|\hfill\\
\hfill\cdot\bigg\|\bigg(\frac{\partial}{\partial\xi}\bigg)^{\beta}(\calG(i\xi))^l\bigg\|_{\Mmm}
\cdot\bigg\|\bigg(\frac{\partial}{\partial\xi}\bigg)^{\alpha}\exp(t\calG(i\xi))\bigg\|_{\Mmm} \,d\xi\\
\hfill\le D_{\varepsilon} e^{(s_0(\calG)+\varepsilon)t}
\sum_{|\alpha|+|\beta|+|\gamma|\le n+2}\frac{(\alpha+\beta+\gamma)!}{\alpha!\beta!\gamma!}
\int_{\Rn}(1+|\xi|^2)^{-j_0-\frac12|\gamma|+ld+k_\alpha}\,d\xi.
\end{multline*}
In the above $d$ is the maximum of the orders of the scalar PDOs which are entries of 
$\calG(\partial_1,\ldots,\partial_n)$. If $j_0$ is sufficiently large, then all the integrals in the last member of the estimate are finite, so that
$(1-\Delta)^{[n/2]+1}({d}/{dt})^l \hat f_t\in L^1(\Rn;\Mmm)$ 
and 
$$
\bigg\|(1-\Delta)^{[n/2]+1}\bigg(\frac{d}{dt}\bigg)^l \hat f_t\bigg\|_{ L^1(\Rn;\Mmm)}
\le K_\varepsilon e^{(s_0(\calG)+\varepsilon)t}
$$
for every $t\in[0,\infty\[$, $l=0,1$ and $\varepsilon>0$,
where $K_\varepsilon\in\]0,\infty\[$ is independent of $t$ and $l$. By \eqref{eq8.1d},
 this implies that $f_t\in L^1(\Rn;\Mmm)$ for every $t\in[0,\infty\[$, and the mapping
$[0,\infty\[\ni t\mapsto f_t\in L^1(\Rn;\Mmm)$ is continuous and satisfies~(c).
\end{proof}

\subsection{Well posedness of the dual spaces}

Let $E$ be an l.c.v.s. continuously imbedded in $\calS^\prime (\Rn;\symC^m)$. Suppose moreover that $\calS (\Rn;\symC^m)$
 is densely and continuously imbedded in~$E$. Let
\begin{multline*}
E^\prime=\Big\{T\in\calS^\prime (\Rn;\symC^m):\mbox{the linear functional}\\
\hfill \calS (\Rn;\symC^m)\ni\varphi\mapsto\<T,\varphi\>
=\sum_{\mu=1}^m T_\mu(\varphi_\mu)\in\symC\\
\hfill\mbox{is continuous on } \calS (\Rn;\symC^m)
\mbox{ in the topology induced by }E\Big\}.
\end{multline*}
Then each $T\in E^\prime$ uniquely extends to a continuous functional on~$E$, and may be identified with that functional. $E^\prime$~is equipped with the topology of uniform convergence on bounded subsets of~$E$.

Assume that $ \calG(\partial_1,\ldots,\partial_n)$ satisfies (ii), and let
$(S_t)_{t\ge0}\subset\calO_C^\prime(\Rn;\Mmm)$ be the 
i.d.c.s. with generating distribution $\calG(\partial_1,\ldots,\partial_n)\otimes\delta$. Then $\calG^{\dag}(-\partial_1,\allowbreak\ldots, -\partial_n)$ also satisfies (ii), 
and $(\check S_t^{\dag})_{t\ge0}\subset\calO_C^\prime (\Rn;\Mmm)$ is the i.d.c.s. with generating distribution $\calG^{\dag}(-\partial_1,\ldots,-\partial_n)\otimes\delta$. ({\spaceskip.33em plus.22em minus.17em The above observation is related to \cite[Sec.~14]{S1}.}) Assume in addition that $\calS (\Rn;\symC^m)$ is dense in the l.c.v.s.~$E$, and 
$E$ is a Montel and hence barrelled space continuously imbedded in $\calS' (\Rn;\symC^m)$. Then
$$
{\rm(iii)}_{\check S^\dag_t,E}\quad\mbox{implies}\quad {\rm(iii)}_{S_t,E^\prime}.
$$
Indeed, the sequential completeness of $E^\prime$  is a consequence of the barrelledness of~$E$. Continuous imbedding of $E^\prime$ in $\calS^\prime (\Rn;\symC^m)$
 follows from dense and continuous imbedding of $\calS (\Rn;\symC^m)$ in~$E$. The other properties of $(S_t)_{t\ge0}$ and $E^\prime$ listed in ${\rm(iii)}_{S_t,E^\prime}$ follow 
from the equality
$$
\<S_t*T,u\>=\<T,\check S_t^\dag * u\>
$$
for all $t\in[0,\infty\[$, $u\in E$ and $T\in E^\prime$.

\subsection{Well posedness of $\calS^\prime (\Rn;\symC^m)$ and uniqueness of solutions}
 If $\calG(\partial_1,\ldots,\partial_n)$ satisfies (ii), then, by the argument presented in Sec.~8.4, $\calS'(\Rn;\symC^m)$ is a well posedness space. Consequently,
$$
((S_t\,*)|_{\calS^\prime (\Rn;\symC^m)})_{t\ge0}\subset
L(\calS^\prime (\Rn;\symC^m);\calS^\prime (\Rn;\symC^m))
$$
is a one-parameter operator semigroup 
of class $(C_{0})$.
Since $\calS^\prime (\Rn;\symC^m)$
is a Montel and hence barrelled space, 
from infinite differentiability of the convolution semigroup $(S_{t})_{t\ge0}\subset \calO^{\prime}_{C}(\Rn;\Mmm)$, it follows by the Banach--Steinhaus theorem that the operator semigroup $((S_{t}\,*)|_{\calS^{\prime}
(\Rn;\Cm)})_{t\ge0}$
is infinitely differentiable in the topology of $L_b(\calS^\prime (\Rn;\symC^m);\allowbreak\calS^\prime (\Rn;\symC^m))$. The infinitesimal generator of this operator semigroup is the differential operator
\begin{equation*}
\calG(\partial_{1},\ldots,\partial_{n})\in L(\calS^\prime (\Rn;\symC^m);\calS^\prime (\Rn;\symC^m)).
\end{equation*}
From the above it is easy to infer that
\begin{align*}
\frac{d}{dt}[(S_t\,*)|_{\calS^\prime (\Rn;\symC^m)}]&= \calG(\partial_1,\ldots,\partial_n)(S_t\,*)|_{\calS^\prime (\Rn;\symC^m)}\\
&= S_t* \calG(\partial_1,\ldots,\partial_n)|_{\calS^\prime (\Rn;\symC^m)}
\end{align*}
for every $t\in[0,\infty\[$, where the derivative is computed in the topology of $L_b(\calS^\prime (\Rn;\symC^m);\calS^\prime (\Rn;\symC^m))$, and is understood as the two-sided derivative if $t\in\]0,\infty\[$, and as the right derivative if $t=0$. {\spaceskip.33em plus.22em minus.17em Consequently, if $u_0\in\calS^\prime (\Rn;\symC^m)$, then the Cauchy problem (iv) has the solution $u(\cdot)\in C^\infty([0,\infty\[;\calS^\prime (\Rn;\symC^m))$ given by formula~(v). If $\tilde u(\cdot)\in C^1([0,\infty\[;\calS^\prime (\Rn;\symC^m))$ is any other solution of (iv), then for every  $t\in\]0,\infty\[$ and $\tau\in\]0,t\[$ one has}
\begin{align*}
\frac d{d\tau}(S_{t-\tau}*\tilde u(\tau))
={}& \lim_{\]0,t-\tau]\ni h\to0}\frac1h(S_{t-\tau-h}-S_{t-\tau})*\tilde u(\tau)\\
&{}+\lim_{\]0,t-\tau]\ni h\to0}S_{t-\tau-h}*\frac d{d\tau}\tilde u(\tau)\\
&{}+ \lim_{\]0,t-\tau]\ni h\to0} S_{t-\tau-h}*
\bigg[\frac 1h(\tilde u(\tau+h)-\tilde u(\tau))-\frac d{d\tau}\tilde u(\tau)\bigg]\\
={}&\bigg(\frac d{d\tau} S_{t-\tau}\bigg)*\tilde u(\tau)+S_{t-\tau}*\frac d{d\tau}
\tilde u(\tau)\\
={}&[S_{t-\tau}*(-\calG(\partial_1,\ldots,\partial_n)\otimes\delta)]
*\tilde u(\tau)\\
&{}+S_{t-\tau}*[(\calG(\partial_1,\ldots,\partial_n)\otimes\delta)*\tilde u(\tau)]=0
\end{align*} 
in the topology of $\calS'(\Rn;\symC^m)$. Indeed, $\calS'(\Rn;\symC^m)$ is barrelled, so that, by the Banach--Steinhaus theorem, the set of convolution operators
$\{S_{t-\tau-h}*:h\in[0,t-\tau]\}$
 is an equicontinuous subset of $L(\calS'(\Rn;\symC^m);\calS'(\Rn;\symC^m))$, whence
$$
\lim_{\]0,t-\tau]\ni h\to 0}S_{t-\tau-h}*\bigg[\frac1h\bigg(\tilde u(\tau+h)-\tilde u(\tau)
-\frac{d}{d\tau}\tilde u(\tau)\bigg)\bigg]=0.
$$
 Consequently, the continuous function
$[0,t]\ni\tau\mapsto S_{t-\tau}*\tilde u(\tau)\in\calS^\prime (\Rn;\symC^m)$ is constant in the open interval $\]0,t\[$, and hence in $[0,t]$. It follows that $\tilde u(t)-S_t* u_0= S_{t-\tau}*\tilde u(\tau)|_{\tau=0}^{\tau=t}=0$. Hence in the class
$C^1([0,\infty\[;\calS^\prime (\Rn;\symC^m))$ the Cauchy problem (iv) with 
$u_0\in\calS^\prime (\Rn;\symC^m)$ has a unique solution, and this solution belongs to $C^\infty([0,\infty\[;\allowbreak\calS^\prime (\Rn;\symC^m))$ and is represented by~(v).

\subsection{Proof of Theorem \ref{thm5}}

Suppose that (ii) holds and $E$ is any l.c.v.s. satisfying ${\rm(iii)}_{S_t,E}$. Then\break
$((S_t\,*)|_{E})_{t\ge0}\subset L(E,E)$ 
is a $(C_0)$-semigroup of operators. Fix $u_0\in E$. The trajectory
$t\mapsto S_t*u_0$ belongs to $C([0,\infty\[;E)$ and, in  view of Sec.~8.5, it
belongs to $C^\infty([0,\infty[;\allowbreak\calS^\prime(\Rn;\symC^m))$ and the right derivative $\frac{d}{dt}\big|_{t=0}(S_t*u_0)$ computed in the topology of
$\calS^\prime(\Rn;\Cm)$ is equal to $\calG(\partial_1,\ldots,\partial_n)u_0$. If the right derivative  $\frac{d}{dt}\big|_{t=0}(S_t*u_0)$ exists in the topology of $E$, then it has the same value $\calG(\partial_1,\ldots,\partial_n)u_0$. This leads to the conclusion that the infinitesimal generator of the operator semigroup 
$((S_t\,*)|_{E})_{t\ge0}\subset L(E;E)$ is equal to the operator $\calG_E$. By 
\cite[Theorem 3.3]{K1} (based on the R.~E.~Edwards boundedness principle for sequentially complete l.c.v.s. \cite[Theorem 7.4.4]{E}), whenever $u_0\in D(\calG_E)$ and
$t\in\]0,\infty\[$, then $S_t*u_0\in D(\calG_E)$ and the two-sided derivative
$\frac d{dt}(S_t*u_0)$ 
computed in the topology of $E$ satisfies the equalities
\begin{equation}\label{eq8.2d}
\frac{d}{dt}(S_t*u_0)=\calG_E(S_t*u_0)=S_t*(\calG_Eu_0).
\end{equation}
Consequently, $u(t)=S_t*u_0$ belongs to $C_1([0,\infty\[;E)$ and is a solution of the Cauchy problem (iv). The uniqueness of this solution is a consequence of Sec.~8.5 and the fact that $u(\cdot)\in C^\infty([0,\infty\[;\calS'(\Rn;\Cm))$.
From \eqref{eq8.2d} it follows that if $u_0\in D(\calG_E^k)$ where  $k=2,3,\ldots,$
then $u(\cdot)\in C^k([0,\infty\[;E)$. 

\section*{Appendix. Proof of (\ref{eq1.28})}

Let  $\calG(\partial_1,\ldots,\partial_n)$ be an $\mm$ matrix whose entries are PDOs
on $\Rn$ with constant complex coefficients. Assume that conditions \eqref{eq1.14} 
and \eqref{eq1.23} are satisfied, so that, by Theorem~\ref{thm4}, 
there is a unique i.d.c.g. $(S_t)_{t\in\symR}\subset \calE'(\Rn;\Mmm)$
with generating distribution $\calG(\partial_1,\ldots,\partial_n)\otimes\delta$. 
Take any $t_0\in\]0,\infty\[$. By the Paley--Wiener--Schwartz theorem, i.e. 
by \cite[Theorem 7.3.1]{H} or \cite[Theorem 8.57]{K-R},
$\scrF S_{t_0}=\exp(t_0\calG(i\,\cdot))$ can be extended to
an $\Mmm$-valued function holomorphic on $\symC^n$ such that for some
$C\in\]0,\infty\[$ and $l\in\symN_0$ one has
$$
\|\!\exp(t_0\calG(i\zeta))\|_{\Mmm}=\|(\scrF S_{t_0})(\zeta)\|_{\Mmm}
\le C(1+|\zeta|)^l e^{H_0(\hIm\zeta)}
$$
for every $\zeta\in\Cn$ where $\Rn\ni \eta\mapsto H_0(\eta)=\sup\{x\eta:
x\in\supp S_{t_0}\}\in\symR$ is the supporting function of
$\supp S_{t_0}$. By \eqref{eq2.1} and \eqref{eq2.2}, it follows that
$$
e^{t_0\max\hRe\sigma(\calG(i\zeta))}\le C(1+|\zeta|)^l e^{H_0(\hIm\zeta)},
$$
and so
\begin{align*}
 \|(\scrF S_{t})(\zeta)\|_{\Mmm}
&=\|\!\exp(t\calG(i\zeta))\|_{\Mmm}\\
&\le(e^{t_0\max\hRe\sigma(\calG(i\zeta))})^{t_0^{-1}t}\bigg(1+\sum_{k=1}^{m-1}\frac{(2t)^k}{k!}
\|\calG(i\zeta)\|_{\Mmm}^k\bigg)\\
&\le C^{t_0^{-1}t}(1\!+\!|\zeta|)^{lt_0^{-1}t}e^{t_0^{-1}t H_0(\hIm\zeta)}
\bigg(\!1\!+\!\sum_{k=1}^{m-1}\!\frac{(2t)^k}{k!}
\|\calG(i\zeta)\|_{\Mmm}^k\bigg)
\end {align*}
for every $(t,\zeta)\in[0,\infty\[\times \Cn$. Since
$\calG(i\zeta)$ is an $\mm$ matrix with polynomial entries, 
it follows that for every $t\in[0,\infty\[$ there are $C_t,k_t\in\]0,\infty\[$ such that
$\|(\scrF S_{t})(\zeta)\|_{\Mmm}\le C_{t}(1+|\zeta|)^{k_t}e^{t_0^{-1}t H_0(\hIm\zeta)}$.
 Hence, by the Paley--Wiener--Schwartz theorem,
$$
\frac 1t\convsupp S_t\subset\frac 1{t_0}\convsupp S_{t_{0}}\ \quad\mbox{for every }
t\in[0,\infty\[.
$$
Since $t$ and $t_0$ can be interchanged, one concludes that
$$
\frac 1t\convsupp S_t=\frac 1{t_0}\convsupp S_{t_{0}}\ \quad\mbox{for every }
t\in[0,\infty\[.\eqno{\rm(A.1)}
$$

The formula
$$
\tilde E(\varphi)=\int_0^\infty S_t(\varphi(t,\cdot))\,dt,\ \quad
\varphi\in\calD(\symR^{1+n}), 
$$
defines a fundamental solution $\tilde E$ for the  matricial PDO \eqref{eq1.26} 
with support contained in the cone
\begin{align*}
K&=\{(t,x_1,\ldots,x_n)\in\symR^{1+n}:t\ge 0,\,(x_1,\ldots,x_n)\in\convsupp S_t\}\\
&=\{(t,x_1,\ldots,x_n)\in\symR^{1+n}:t\ge 0,\,(x_1,\ldots,x_n)\in tt_{0}^{-1}\convsupp S_{t_0}\}
\end{align*}
where the equality is a consequence of~(A.1). 
Consequently, $E=\det_*\tilde E$ is a fundamental solution for the operator
$P(\partial_1,\partial_1,\ldots,\partial_n)=\det(\symjed_{\mm}\otimes
\partial_t-\calG(\partial_1,\ldots,\partial_n))$. 
In the above $\det_*$ is the determinant in the sense of the convolution algebra $\calO_C^\prime(\Rn)$. It follows that
$\supp E\subset K$. Since $K\subset H_+:=\{(t,x_1,\ldots,x_n)\in
\symR^{1+n}:t\ge 0\}$, from \cite[Theorem 12.5.1]{H} it follows that
\begin{equation}
\supp E\subset \varGamma^0,\tag{A.2}
\end{equation}
and
%\Anumer
$$ 
\textbot{whenever $H$ is a convex cone with
$\supp E\subset H\subset H_+,$ then $\varGamma^0\subset H.$}
\eqno{\rm(A.3)}
$$
%\enumer
From (A.3) it follows that
$$
\varGamma^0\subset K.
\eqno{\rm(A.4)}
$$
		In order to prove the inclusion opposite to (A.4), first we shall show that
$$
\supp \tilde E\subset\varGamma^0. \eqno{\rm(A.5)}
$$
The above inclusion follows from (A.2) and the equality
$$
\tilde E=[\adj(\symjed_{\mm}\otimes \partial_t-\calG(\partial_1,\ldots,\partial_n))]
\otimes E.
\eqno{\rm(A.6)}
$$

Indeed,
$E_1=\tilde E$ and $E_2=[\adj(\symjed_{\mm}\otimes\partial_t
-\calG(\partial_1,\ldots,\partial_n))]\otimes E$
both have support contained in $H_+$, and both are fundamental solutions for the 
matricial PDO \eqref{eq1.26}. Moreover $\vartheta_0 E_i\in\calE'(\symR^{1+n};\Mmm)$ for
$i=1,2$ and every $\vartheta_0\in C^\infty(\symR^{1+n})$
 such that $\vartheta_{0}(t,x_1,\ldots,x_n)\equiv \vartheta(t)$ where $\vartheta\in \calD(\symR)$.
 These properties of  $E_i$, $i=1,2$, imply the equality $E_1=E_2$ (see the
author's preprint {\it The Petrovski\u\i\ condition and rapidly decreasing distributions}, Inst. Math., Polish Acad. Sci., 2011). 
The equality $E_1=E_2$ means that (A.6) holds. Now, (A.5) is a consequence of (A.2) and (A.6).

From (A.5) the inclusion
$$
K\subset \varGamma^0\eqno{\rm(A.7)}
$$		
may be deduced by an elementary reasoning. Indeed, (A.7) follows once it is proved that
$$
\supp S_t\subset \varGamma_t^0:=\{(x_1,\ldots,x_n)\in\Rn:
(t,x_1,\ldots,x_n)\in\varGamma^0\}\eqno{\rm(A.8)}
$$		 
for every $t\in[0,\infty\[$. If (A.8) were not true for some $t_0\in[0,\infty\[$, 
then there would exist $\varphi_0\in \calD(\Rn)$
such that $(\supp\varphi_0)\cap \varGamma_{t_0}^0=\emptyset$ and
$0\ne S_{t_0}(\varphi_0)\in\Mmm$. Since
$S_t(\varphi_0)$ depends continuously on $t$, there would exist
$\psi_0\in \calD(\]0,\infty\[)$ such that
$(\supp(\psi_0\otimes \varphi_0))\cap \varGamma^0=\emptyset$
and $\tilde E(\psi_0\otimes \varphi_0)=\int_0^\infty\psi_0(t)S_t(\varphi_0)\,dt\ne0$, 
contrary to (A.5). Therefore (A.8) is true, and (A.7) holds. The inclusions (A.4) and (A.7) prove the equality (1.28).

\end{document}